%% file: autorank2.tex
\documentstyle[epsf,12pt]{article}
\newtheorem{proposition}{~~~~Proposition}[section]

\newtheorem{lemma}{~~~~Lemma}
\newtheorem{corollary}{~~~~Corollary}
\newtheorem{theorem}{~~~~Theorem}
\newtheorem{definition}{~~~~~~{\sc Definition}}
\newtheorem{remark}{~~~~~~{\sc Remark}}
\newcommand{\hako}[1]{\mbox{\boldmath$#1$}}
\def\spmapright#1{\smash{%$w_k^{(j)} =i$
 \mathop{\hbox to 1.2cm{\rightarrowfill}}
  \limits^{#1}}}
\def\sbmapright#1{\smash{%
 \mathop{\hbox to 1.2cm{\rightarrowfill}}
  \limits_{#1}}}
\def\lmapdown#1{\Big\downarrow
 \llap{$\vcenter{\hbox{$\scriptstyle#1\,$}}$ }}   
\def\rmapdown#1{\Big\downarrow
 \rlap{$\vcenter{\hbox{$\scriptstyle#1$}} $}}  

\title{ Stepped surfaces and Rauzy fractals induced from automorphisms
on the free group of rank 2}
\author{{\sc Hiromi} EI }
\date{}
\begin{document}
%\setcounter{page}{0}
%~~~
%\clearpage
\maketitle
\thispagestyle{myheadings}
\setcounter{section}{-1}
%\markright{Stepped surfaces and Rauzy fractals induced from automorphisms}
%on the free group of rank 2 in some class
%\keywords{Keyword one, keyword two, keyword three}
\begin{abstract}
For substitution satisfying Pisot, irreducible, unimodular condition,
a tiling substitution plays a key role in construction 
of a stepped surface and Rauzy fractal (see \cite{a-i}).
In this paper we will extend the method to hyperbolic automorphisms on the 
free group of rank 2 in some class, and obtain set equations of Rauzy
 fractals by virtue of a tiling substitution. We will also see that 
the domain exchange transformation on Rauzy fractal 
is just a two interval exchange transformation. 
\end{abstract}

\def\theth1{\relax}
\input{autorank2-0.tex}
\input{autorank2-1.tex}
\input{autorank2-2.tex}
\input{autorank2-3.tex}
\input{autorank2-4.tex}
\\
%\cite{a-i} \cite{ak2}.
%\input{ref.tex}
%\bibliographystyle{abbrvnat}
% use the following instead if you encounter problems 
\bibliographystyle{plain}
\nocite{i-k}
\nocite{i-r-e}
\nocite{b-f}
\nocite{a-b-h-s}
\nocite{i-y}
\nocite{i-f-h}
\nocite{f}
\nocite{en}
\bibliography{autorank2}
% to list all of references inref.bib 
% \nocite{*}
\label{sec:biblio}

\end{document}

%% file: autorank2-0.tex
\section{Introduction}

Rauzy fractal \cite{rauzy} has been extensively studied because it
plays significant roles in the study of substitutive dynamical system
in the case of Pisot, irreducible unimodular substitution
(e.g., \cite{a-i, i-r, b-f-m-s}).
Arnoux and Ito \cite{a-i} gives the way to construct a stepped surface
(see Proposition \ref{prop:0-1} and Figure \ref{fig:1}) 
and Rauzy fractal (see Proposition \ref{prop:0-2} and
Figure \ref{fig:2-0})  by using a tiling substitution 
which is sometimes called 
a dual map, and the set equation of Rauzy fractal related to its
self-similality (see Proposition \ref{prop:0-3}); and they
obtains the domain exchange transformation on Rauzy fractal as the 
realization of the substitutive dynamical system 
(see Theorem \ref{theorem:0-1}).
As the extension to automorphisms on the free group, 
Arnoux, Berht{\'e}, Hilion and 
Siegel \cite{a-b-h-s} have started the study of the class 
where cancellation of letters under the iteration of automorphism 
does not occur; 
and Berht{\'e} and Fernique \cite{b-f} discussed the action of a tiling
substitution for automorphism on the stepped surface.

In this paper we study the natural class of automorphisms related to the
companion
matrices of quadratic polynomials $x^2 - ax \mp 1$ such that 
\begin{equation}
\label{eq:1}
A_{\pm}=\pmatrix{ 0 & \pm 1  \cr 1 & a},
\end{equation}
and assume ``hyperbolicity'' instead of the Pisot condition.
The purpose of this paper is to find automorphisms related to the
matrices given by (\ref{eq:1}) which give analogue properties in 
substitutions case, 
and to discuss stepped surfaces, Rauzy fractals and dynamical systems
on the chosen automorphisms.
So we shall extend the way and technique for substitutions to 
ones for automorphisms. But we sometimes encounter the problem which is 
peculiar to automorphisms.
For example, take the automorphism defined by 
\[
\sigma :~
\left\{\
\begin{array}{l}
1 \to 2\\
2 \to 21^{-1}22
\end{array}
\right. ,
\]
then 
cancellation of letters under the iteration of $\sigma$ occurs
because
\[
 \sigma^2(2)=\sigma(21^{-1}22)=21^{-1}22 ~ \underline{2^{-1} ~ 2}1^{-1}22 ~ 21^{-1}22.
\]
Such cancellation never occurs for any substitution and automorphisms
discussed in \cite{a-b-h-s}.
Main idea to solve this problem is to find 
a substitution or a ``pseud-substitution'' $\tau$, which is called 
an ``alternative substitution'' in this paper, for each
automorphism $\sigma$ satisfying
$\sigma = \delta ^{-1} \circ \tau \circ \delta$ with some automorphism
$\delta$.
We will show that many results obtained in substitution case
also hold for the chosen automorphisms by using conjugate $\tau$.
%to which the automorphism is conjugate.

In Section 1 recalls results in the case of substitutions of rank 2.
So similar results will appear
in the case of some automorphisms. 

Under the condition of hyperbolicity, there are four cases 
of the matrix given by (\ref{eq:1}).
In Section 2, we choose automorphisms on the free group of rank 2
for each cases, and find their conjugates which are substitutions or 
alternative substitutions.
These automorphisms are discussed in the following sections.

In Section 3, 
we show that stepped surfaces related to the
chosen automorphisms can be obtain by ones related to their conjugates.
By using this fact, we will find 
appropriate initial elements, so called seeds, 
for tiling substitutions, 
and generate the stepped surfaces related to the automorphisms.

In Section 4 is devoted to Rauzy fractals induced from the automorphisms. 
First we generate Rauzy fractals for 
the both of the automorphisms and its conjugates
by each tiling substitutions with appropriate seeds;
and show that Rauzy fractals induced from the automorphisms 
can be written as a disjoint union of Rauzy fractals 
related to their conjugates, and thus they are just intervals
in Theorem \ref{theorem:3-1}. 
Second we consider measurable dynamical systems 
with domain exchange transformations on 
Rauzy fractals, and the structure of its induced transformations
in Theorem \ref{theorem:3-2}. 
Finally we see the Rauzy fractals related to the automorphisms
are obtained by their fixed points or periodic points
in Theorem \ref{theorem:3-3}. .
%In Section 1, for chosen automorphisms 
%we find conjugacies. 
%In Section 2, we generate the stepped surface with
%a tiling substitutions. Here we also know the way to determin a seed
%to generate a stepped surface.
%In Section 3, we discuss the symbolical dynamical system generated
%by a fixed point or periodic point of an automorphism and the domain
%exchange transformation.

%% file: autorank2-1.tex
\section{Results in substitution case}

We briefly recall the substitution case.
We concentrate substitutions of rank 2 even though some 
properties are true for any rank.
%, where rank means the number of letters in ${\cal A}$. 
Let ${\cal A}=\{ 1,2\}$ (resp. $\widehat{\cal A} =\{ 1,2,1^{-1},2^{-1}\}$ )
be an alphabet consisting of two letters (resp. four letters), 
and ${\cal A}^*$ (resp. $\widehat{\cal A}^*$) the free monoid
with the empty word $\epsilon$ generated by ${\cal A}$
(resp. $\widehat{\cal A}$ ).
More preciously, a word $W=w_1 w_2 \cdots w_n  \in {\cal A}^*$ (resp. 
$\widehat{\cal A}^*$)
satisfies $w_i \in {\cal A} $ 
(resp. $w_i \in \widehat{\cal A}$ ) for any $i \in \{1,2,\cdots , n\}$.
We say a word $W=w_1 w_2 \cdots w_n \in \widehat{\cal A}^*$ is reduced if 
$w_i w_{i+1} \ne\epsilon$ for any $i \in \{1,2,\cdots , n-1\}$. 
A word $w_1 w w^{-1} w_2  \in  \widehat{\cal A}^*$ becomes $w_1 w_2$ after
cancellation. If two words $W_1,~W_2 \in \widehat{\cal A}^*$  becomes
the same reduced word after cancellation, then we say they are referred to be
equivalent, and written as $W_1 \sim W_2$. 
The free group of rank 2 is defined by $F_2= \widehat{\cal A}^* / \sim$.
For simplicity, the concatenation of $k$ copies of some letter $i \in
{\cal A}$
(resp. $i^{-1} \in \{1^{-1},2^{-1} \}$)
is written as $ii\cdots i=i^k$ (resp. $i^{-1}i^{-1}\cdots i^{-1}=i^{-k}$).
An endomorphism on 
${\cal A}^*$ is called a {\it substitution} of rank 2 over ${\cal A}$ and it is
naturally extended to an endomorphism on $F_2$. 
A substitution is referred to be {\it invertible} if it is an automorphism on
$F_2$ by the extension. 
%From now on, to avoid confusion, we use the notation $\tau$
%(resp. $\sigma$) 
%for a substitution or a alternative substitution (resp. an endomorphism on
%the free group $F_2$).

A canonical homomorphism ${\bf f}: 
F_2 \to {\mbox{\boldmath $Z$}}^2$ 
is defined by ${\bf f}(\epsilon)={\hako o}$ and
${\bf f}(i^{\pm 1})=\pm {\hako{e}}_i,~i \in {\cal A}$. 
Then for a matrix $A_{\sigma}$ defined by 
$({\bf f}(\sigma (1)),{\bf f}(\sigma (2)))$, 
so called an {\it incidence matrix} associated with an endomorphism $\sigma$
on $F_2$,
 the following diagram becomes commutative:
\[
\begin{array}{ccc}
 F_2 & \spmapright{\sigma} & F_2  \\
\lmapdown{\bf f} &  &\rmapdown{\bf f}\\
 {\mbox{\boldmath $Z$}}^2 & \sbmapright{A_{\sigma} }& 
 {\mbox{\boldmath $Z$}}^2
\end{array}~~~.
\]

\begin{figure}[hbtp]
\setlength{\unitlength}{1mm}
\begin{center}
\begin{picture}(40,15)(0,0)
\put(10,5){\line(0,1){10}}
\put(0,15){\line(1,0){10}}
\put(0,5){\circle*{1}}
\put(0,0){${\cal U}$}
\put(20,5){\line(0,1){10}}
\put(20,5){\line(1,0){10}}
\put(20,5){\circle*{1}}
\put(20,0){${\cal U}'$}
\end{picture}
\epsfxsize=5cm
\epsfbox{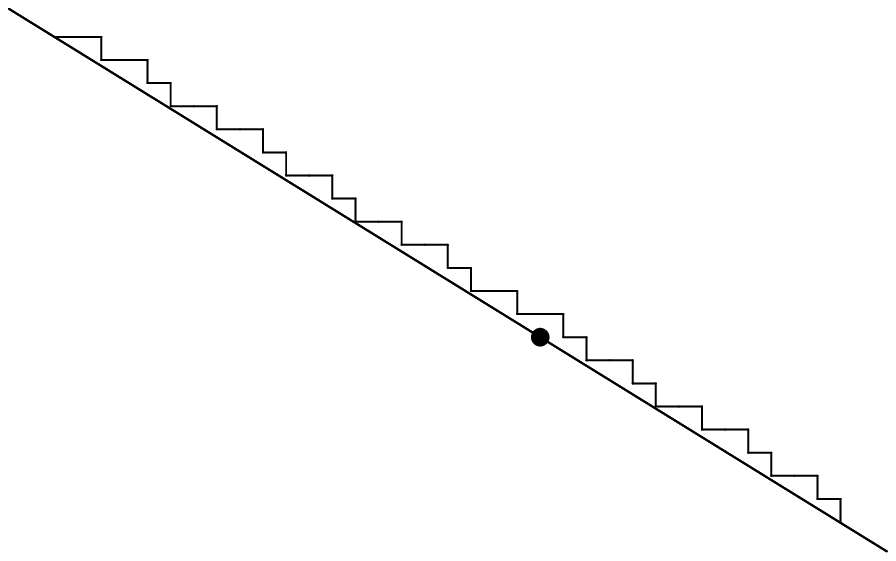}
\epsfxsize=5cm
\epsfbox{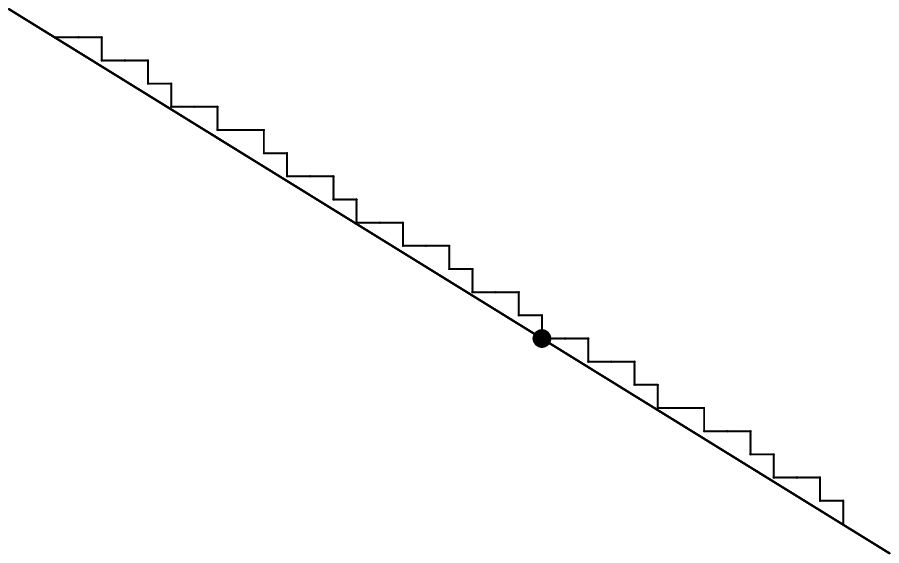}
\end{center}
\caption{The seeds ${\cal U},~{\cal U}' $ and $\sigma^{*~n}({\cal
 U}),~\sigma^{*~n}({\cal U}')$ for the substitution given by (\ref{sub})}\label{fig:1}
\end{figure}

For example, let us consider the substitution $\sigma$ of rank 2 given by
\begin{eqnarray}\label{sub}
\sigma :~
\left\{\
\begin{array}{l}
1 \to 2\\
2 \to 21
\end{array}
\right. ,
\end{eqnarray}
with the incidence matrix 
\[
A_{\sigma}=\pmatrix{ 0 & 1  \cr 1 & 1}.
\]
The substitution is Pisot, irreducible, unimodular,
that is, the characteristic polynomial $\Phi_{\sigma}(x)$ of
$A_{\sigma}$ satisfies the following three conditions:
\begin{itemize}
\item (Pisot condition)~The maximum root of $\Phi_{\sigma}(x)$ is Pisot
 number, that is, the dominant eigenvalue of $A_{\sigma}$ is greater
than one and the other has modulus less than one,
\item (Irreducible condition)~$\Phi_{\sigma}(x)$ is irreducible over
      ${\mbox{\boldmath $Q$}}$,
\item (Unimodular condition)~$\mid \det A_{\sigma} \mid =1$.
\end{itemize}
A substitution $\sigma$ is referred to be {\it primitive} if there exists $n$ such that
for any pair $(i,j)$ the letter $i$ occurs in the words $\sigma^n(j)$, in
other words, the incidence matrix $A_{\sigma}$ of $\sigma$ is primitive.
In this section we assume a substitution is Pisot, irreducible,
unimodular and primitive. 
By Perron-Frobenius Theorem and the Pisot condition,
the incidence matrix $A_{\sigma}$ of 
a Pisot irreducible, unimodular and primitive substitution $\sigma$
has a positive column eigenvector ${\hako u}_{\sigma} >0$,
a positive lower eigenvector ${\hako v}_{\sigma} >0$
corresponding to the positive eigenvalue $\lambda_{\sigma} >1$,
and another column eigenvector ${\hako u}_{\sigma}'$ 
corresponding to the other eigenvalue 
$\lambda_\sigma ' $ with $|\lambda_\sigma ' |<1$.
It is easy to check the contractive eigenspace $P_{\sigma}$ of
$A_{\sigma}$ spanned by ${\hako u}_{\sigma}'$ is given by $P_{\sigma}=
\{ {\hako x}\in {\mbox{\boldmath $R$}}^2 \mid <{\hako x},
{}^t{\hako v}_{\sigma}>=0\}$, where $<\cdot,\cdot>$ is an inner product; 
and the stepped surfaces of $P_\sigma$ are defined by
\begin{eqnarray*}
 {\cal S}_{\sigma}&:=& \bigcup_{({\hako x}, i^*) \in S_{\sigma}} ({\hako x}, i^*),\\
{\cal S}_{\sigma}'&:=& \bigcup_{({\hako x}, i^*) \in S_{\sigma}'} ({\hako x}, i^*),
\end{eqnarray*}
where
\begin{eqnarray*}
 S_{\sigma} &:=& \left\{({\hako x}, i^*) \in  {\mbox{\boldmath $Z$}}^2 \times \{1^*,2^*\}  
\mid 
\left<{\hako x},{}^t{\hako v}_{\sigma} \right> > 0,~  
\left<{\hako x}-{\hako e}_i ,{}^t{\hako v}_{\sigma}\right> \le 0\right\},\\
 S_{\sigma}'&:=& \left\{({\hako x}, i^*) \in  {\mbox{\boldmath $Z$}}^2 \times \{1^*,2^*\}  
\mid 
\left<{\hako x},{}^t{\hako v}_{\sigma}\right> \ge 0,~
\left<{\hako x}-{\hako e}_i ,{}^t{\hako v}_{\sigma}\right> < 0\right\}.
\end{eqnarray*}
Identify $({\hako x}, i^*) \in {\mbox{\boldmath $Z$}}^2 \times \{1^*,2^*\}$
with the positive oriented unit segment spanned by the 
fundamental vector ${\hako e}_j$ translated by ${\hako x}$, 
where $\{i,j\}=\{1,2\}$ (see Figure \ref{fig:3}), then
these stepped surfaces ${\cal S}_{\sigma}, {\cal S}_{\sigma}'$ are 
discrete approximations of $P_{\sigma}$ (see Figure \ref{fig:1}).
On the other hand, the stepped surface is generated by
a tiling substitution.
On the free ${\mbox{\boldmath $Z$}}$-module ${\cal G}^*$ defined by
\[
{\cal G}^*:=\left\{ \sum_{k=1}^l n_k ({\hako x}_k, i_k^*) \mid 
n_k \in {\mbox{\boldmath $Z$}},~ 
{\hako x}_k \in {\mbox{\boldmath $Z$}}^2,~
i_k \in {\cal A} \mbox{ for any } k,~l < \infty \right\}, 
\]
an endomorphism $\sigma^*$, so called a
tiling substitution, is given by
\[
\sigma^*({\hako x} , i^*) = \sum_{j \in {\cal A}} \sum_{w_k^{(j)}=i} 
(A_{\sigma}^{-1} ({\hako x}+{\bf f}(S_k^{(j)})) , j^*)~\mbox{for~}
({\hako x} , i^*) \in {\cal G}^* ,
\]
where ${\bf f}$ is a canonical homomorphism from the free monoid 
${\cal A}^*$ to $ {\mbox{\boldmath $Z$}}^2$, $\sigma(j)=w_1^{(j)}
w_2^{(j)} \cdots w_{l^{(j)}}^{(j)}$ and $S_k^{(j)}=w_{k+1}^{(j)}
w_{k+2}^{(j)} \cdots w_{l^{(j)}}^{(j)}$.
%When we identifies a substitution with an endomorphism on the free
%group, the given substitution is invertible.
Remark that we usually use the notations ${\cal G}_1^*,~E_1^*(\sigma)$
instead of ${\cal G}^*,~\sigma^*$ when we consider substitutions of higher rank
(cf. \cite{a-i, s-a-i}).
For the substitution $\sigma$ given by (\ref{sub}),
the tiling substitution $\sigma^*$ is determined as follows:
\[
 \sigma^* ({\hako x}, i^*)=
\left\{
\begin{array}{ll}
(A_{\sigma}^{-1} {\hako x} , 2^*) & \mbox{if } i=1\\
(A_{\sigma}^{-1} {\hako x} , 1^*)
+(A_{\sigma}^{-1}{ \hako x} - {\hako e}_1+ {\hako e}_2, 2^*)& \mbox{if } i=2\\
\end{array}
\right. .
\]
We also identify an element of ${\cal G}^*$
with a union of oriented unit segments with multiplicity.
Define the subset of ${\cal G}^*$ which consists of unit segments on
the stepped surface without multiplicity as follows:
\[
 {\cal G}_{\sigma}^*:=
\left\{
\sum_{k=1}^l({\hako x}_k, i_k^*) \mid 
\begin{array}{l}
({\hako x}_k, i_k^*) \in {S_{\sigma}},~l<\infty\\
({\hako x}_k, i_k^*) \ne ({\hako x}_{k'}, i_{k'}^*) \mbox{ if } k \ne k'
\end{array}
\right\},
\]
and $ {{\cal G}_{\sigma}^*}'$ is also defined in the same way by replacing 
${S_{\sigma}}$ with $S_{\sigma}'$.
By iterating $\sigma^*$ for the initial elements 
${\cal U}:=({\hako e}_1,1^*)+({\hako e}_2, 2^*)\in {\cal G}_{\sigma}^*,~
{\cal U}':=({\hako o},1^*)+({\hako o}, 2^*) 
\in {{\cal G}_{\sigma}^*}'$, the stepped surfaces are obtained.
\begin{proposition}
\label{prop:0-1}
(\cite{a-i})
For a substitution $\sigma$,
we have 
$\sigma^{*~n} ({\cal U}) \in {\cal G}_{\sigma}^*$ (resp. $\sigma^{*~n}
 ({\cal U}')  \in {{\cal G}_{\sigma}^*}'$)
and $\sigma^{*~n}({\cal U}) - \sigma^{*~n}
 ({\cal U}')= {\cal U}-{\cal U}'$ for any positive integer $n$.
%, and
%there is no overlap in $E_1^{*~n}(\sigma) ({\cal U}),~E_1^{*~n}(\sigma)
% ({\cal U}')$ in the sense of Lebesgue measure. 
\end{proposition}
By using the projection $\pi_{\sigma}$ from ${\mbox{\boldmath $R$}}^2$ to 
$P_{\sigma}$ along the eigenvector ${\hako u}_{\sigma}$, 
we obtain the quasi-periodic tiling ${\cal T}$ on $P_{\sigma}$ with
two prototiles:
\[
 {\cal T} :=  \bigcup_{({\hako x}, i^*) \in S_{\sigma}}\pi_{\sigma} ({\hako x}, i^*).
\]
That is the reason why $\sigma^*$ is called a tiling substitution.

We call ${\cal U},~{\cal U}'$ ``seeds'' for the tiling substitution
$\sigma^{*}$. The choice of seeds is important when we consider 
Rauzy fractals and  dynamical systems on them generated by a
substitution. 
Recall that we identify an element 
$\sum_{k=1}^l ({\hako x}_k, i_k^*) \in {\cal G}_{\sigma}^*$ with 
$\cup_{k=1}^l ({\hako x}_k, i_k^*) \subset {\cal S}_{\sigma}$.
\begin{proposition}
\label{prop:0-2}
(\cite{a-i})
There exist the following limit sets in the sense of Hausdorff metric
for $i \in {\cal A}$:
\begin{eqnarray*}
X_{\sigma}&:=& \lim_{n \to \infty} A_{\sigma}^n \pi_{\sigma}
 \sigma^{*~n}({\cal U}),\\
&=& \lim_{n \to \infty} A_{\sigma}^n \pi_{\sigma}
 \sigma^{*~n}({\cal U}'),\\
X_{\sigma}^{(i)}&:=& \lim_{n \to \infty} A_{\sigma}^n \pi_{\sigma}
 \sigma^{*~n}({\hako e}_i, i^*),\\
{X'}_{\sigma}^{(i)}  &:=& \lim_{n \to \infty} A_{\sigma}^n \pi_{\sigma}
 \sigma^{*~n}({\hako o}, i^*).
\end{eqnarray*}
Since the boundaries 
of the sets $X_{\sigma},~X_{\sigma}^{(i)},~{X'}_{\sigma}^{(i)}$
are fractal in the case of substitutions of higher rank,
so they are called Rauzy fractals or atomic surfaces.
\end{proposition}
%%%%%%%
\begin{figure}[hbtp]
\begin{center}
\begin{minipage}{6cm}
~~~~~~~~$X_\sigma^{(1)}$~~~~~~~~~~$X_\sigma^{(2)}$\\
\epsfxsize=5cm
\epsfbox{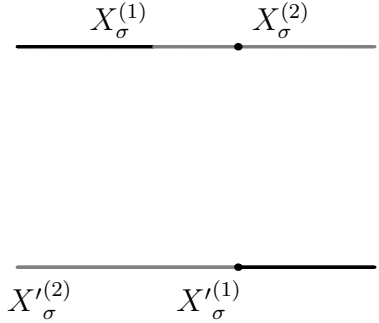}
\\
~~~~~~~~${X'} _\sigma^{(2)}$~~~~~~~~~~${X'} _\sigma^{(1)}$
\end{minipage}
\end{center}
\caption{Rauzy fractals $X_{\sigma}^{(i)},~{X'} _{\sigma}^{(i)},~i \in
 {\cal A}$ related to the substitution $\sigma$ given by
 (\ref{sub})}\label{fig:2-0}
\end{figure}
It is well known that these Rauzy fractals are given by a fixed point
of a substitution as follows:
\begin{eqnarray}
\label{atomic1}
X_{\sigma}^{(i)}&=& 
\overline{ \{-\pi_{\sigma}{\bf f} (s_0 s_1 \cdots s_{k-1}) | s_k=i\}},\\
\label{atomic2}
{X'}_{\sigma}^{(i)}&=&
\overline{ \{-\pi_{\sigma}{\bf f} (s_0 s_1 \cdots s_{k}) | s_k=i\}},
\end{eqnarray}
where the one-sided sequence $s_0 s_1 \cdots$ is a fixed point or a
periodic point of a substitution $\sigma$ and $\overline{A}$ means the closure of $A$
(cf. \cite{a-i, i-r}).

By the definition of $\sigma^*$ and
$X_{\sigma}^{(i)},~{X'}_{\sigma}^{(i)}$,
we have the proposition:
\begin{proposition}
\label{prop:0-3}
(\cite{a-i})
The following set equations hold for $i \in {\cal A}$:
\begin{eqnarray*}
A_{\sigma}^{-1}X_{\sigma}^{(i)}  &= & \cup_{j \in {\cal A}}
 \cup_{w_k^{(j)}=i} 
(-A_{\sigma}^{-1}\pi_{\sigma} {\bf f}(P_k^{(j)}) + X_{\sigma}^{(j)} ),\\
A_{\sigma}^{-1}{X'}_{\sigma}^{(i)}  &= &
\cup_{j \in {\cal A}} \cup_{w_k^{(j)}=i} 
(A_{\sigma}^{-1} \pi_{\sigma}{\bf f}(S_k^{(j)}) + {X'}_{\sigma}^{(j)} ),
\end{eqnarray*}
where ${\hako x}+S:=\{{\hako x}+{\hako y}|~ {\hako y} \in S\}$ for $S
 \subset P_\sigma,~{\hako x} \in P_\sigma$.

Moreover, the sets $(-A_{\sigma}^{-1}  \pi_{\sigma}{\bf f}(P_k^{(j)}) +
 X_{\sigma}^{(j)}), ~{j \in {\cal A}}$ such that 
$w_k^{(j)}=i$ are disjoint in the sense of Lebesgue
measure, and the same holds true for the sets $(-A_{\sigma}^{-1}  
\pi_{\sigma}{\bf f}(S_k^{(j)}) + {X'}_{\sigma}^{(j)})$.
\end{proposition}

\begin{figure}[hbtp]
\begin{center}
\begin{minipage}{6cm}
~~~$X_\sigma^{(1)}$~~~~$X_\sigma^{(2)}$\\
\epsfxsize=5cm
\epsfbox{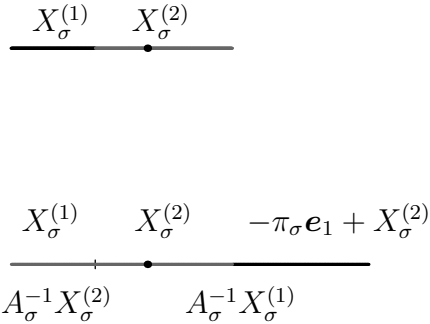}
\\
~~~~~~$A_\sigma^{-1} X_\sigma^{(2)}$~~~~~~~$A_\sigma^{-1} X_\sigma^{(1)}$
\end{minipage}
\setlength{\unitlength}{1mm}
\begin{picture}(1,1)(0,0)
\put(-60,-8){$X_\sigma^{(1)}$}
\put(-45,-8){$X_\sigma^{(2)}$}
\put(-30,-8){$-\pi_\sigma {\hako e}_1 + X_\sigma^{(2)}$}
\end{picture}
\end{center}
\caption{The set equations
 $A_{\sigma}^{-1}X_{\sigma}^{(1)}=-\pi_\sigma
{\hako e}_1 +X_{\sigma}^{(2)},~
A_{\sigma}^{-1}X_{\sigma}^{(2)}=X_{\sigma}^{(1)} \cup  X_\sigma^{(2)}$
 for the substitution given by (\ref{sub})}\label{fig:2}
\end{figure}
\begin{figure}[hbtp]
\begin{center}
\begin{minipage}{6cm}
~~~~~~~~~~~~~~$A_{\sigma}X_\sigma^{(2)}$~~~~~$A_{\sigma}X_\sigma^{(1)}$\\
\epsfxsize=5cm
\epsfbox{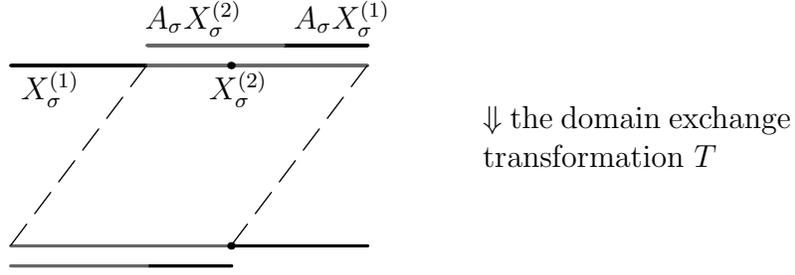}
%\\
%~~~~~~~~${X'} _\sigma^{(2)}$~~~~~~~~~~${X'} _\sigma^{(1)}$
\end{minipage}
\setlength{\unitlength}{1mm}
\begin{picture}(0,0)(0,0)
\put(-60,6){$X_\sigma^{(1)}$}
\put(-35,6){$X_\sigma^{(2)}$}
\end{picture}
\begin{minipage}{4.1cm}
$\Downarrow$ the domain exchange\\
~~~~ transformation $T$
\end{minipage}
\end{center}
\caption{The domain exchange transformation $T$ for the substitution given by (\ref{sub})}\label{fig:3-2}
\end{figure}

\begin{definition}
Let $(X,T,\mu)$ be a measurable dynamical system, $\sigma$ a substitution over 
the alphabet ${\cal A}$  such that
\[ \sigma (i) = w_1^{(i)}w_2^{(i)}\cdots w_{l^{(i)}}^{(i)} ~, \]
$\{ X^{(i)} \mid i \in {\cal A} \}$ a measurable partition  of $X$,
and $\{ A^{(i)} \mid i \in {\cal A} \}$ 
a measurable partition of a subset $A$ of $X$.
We say that the transformation $T$ has {\bf $\sigma$-structure} with respect to
the pair of partitions $\{ X^{(i)} \},~\{ A^{(i)} \}$ if 
the following conditions hold up to set of measure $0$:
\begin{center}
\(
\begin{array}{ll}
 T^{k} A^{(i)}  \subset X^{(w_{k+1}^{(i)})}
&~for~all~~ i\in {\cal A}~,~k=0,1, \cdots, l^{(i)}-1\\
T^kA^{(i)}  \cap A = \emptyset
 &~for~all~~i \in {\cal A},~ 0<k<l^{(i)}\\
T^{l^{(i)}} A^{(i)} \subset A &~for~all~i\in {\cal A}\\
X=\bigcup_{i \in {\cal A}} \bigcup_{0 \leq k \leq l^{(i)}-1 } T^k A^{(i)}
&~(non-overlapping)
\end{array}.
\)
\end{center}
\end{definition}
\begin{theorem}
\label{theorem:0-1}(\cite{a-i})
For a Pisot, unimodular, irreducible and primitive substitution, 
define the map $T:~ X_{\sigma} \to X_{\sigma}$ by
\[
T({\hako x}):= {\hako x} - \pi_{\sigma}({\hako e}_i)~\mbox{if } {\hako
 x} \in  X_{\sigma}^{(i)}.
\]
The map $T$, so called a domain exchange transformation, is
 well-defined; and the measurable dynamical system $(X_{\sigma},T,\mu)$
 with Lebesgue measure $\mu$ has $\sigma^n$-structure with respect to the pair
 of partitions $\{X_{\sigma}^{(i)}|~i \in {\cal A}\}$, 
$\{A_{\sigma}^n X_{\sigma}^{(i)}|~i \in {\cal A}\}$ (see Figure \ref{fig:3-2}).
\end{theorem}
For the substitution given by (\ref{sub}), 
since $A_\sigma X_\sigma^{(1)} \subset X_\sigma^{(2)},~ 
T(A_\sigma X_\sigma^{(1)}) \subset A_\sigma X_\sigma$ and
$A_\sigma X_\sigma^{(2)} \subset X_\sigma^{(2)},~ 
T(A_\sigma X_\sigma^{(2)}) \subset X_\sigma^{(1)},~
T^2(A_\sigma X_\sigma^{(2)}) \subset A_\sigma X_\sigma$,
we can check the domain exchange transformation $T$ has
$\sigma$-structure, and moreover $\sigma^n$-structure
with respect to the partitions 
 $\{X_{\sigma}^{(i)}|~i \in {\cal A} \}$, $\{A_{\sigma}^n 
X_{\sigma}^{(i)}|~i \in {\cal A} \}$ for any positive integer $n$.
On the other hand, the fixed point $s_0s_1\cdots = 
\lim_{n \to \infty} \sigma^n(2)$ and the origin ${\hako o} 
\in A_{\sigma}^n X_{\sigma}^{(2)}$ for any positive integer $n$. 
Therefore from $\sigma^n$-structure, we have 
\[
 T^k({\hako o}) \in X_{\sigma}^{(s_k)}\mbox{ for all }k=0,1, \cdots.
\]

At the end of review of the case of substitutions, recall the following theorem 
related to the topological property of Rauzy fractals.
\begin{theorem}
\label{theorem:0-2}
(\cite{b-e-i-r})
Let a substitution $\sigma$ of rank 2 be  Pisot, unimodular, irreducible and
 primitive.
The Rauzy fractals $X_{\sigma},~X_{\sigma}^{(i)},~{X'}_{\sigma}^{(i)},~
i \in {\cal A} $
are interval 
if and only if $\sigma$ is invertible. Moreover, 
$X_{\sigma}^{(i)},~{X'}_{\sigma}^{(i)}$ are intervals given by
\begin{eqnarray*}
X_{\sigma}^{(i)}&=&\pi_{\sigma} ({\hako e}_i,i^*) + {\hako h}\\
{X'}_{\sigma}^{(i)}&=&\pi_{\sigma} ({\hako o},i^*) + {\hako h}
\end{eqnarray*}
for some ${\hako h} \in P_{\sigma}$.
\end{theorem}

From the theorem, if $\sigma$ is invertible, then the domain exchange
transformation $T$ is just a two interval exchange transformation on the
one dimensional torus.

%% file: autorank2-2.tex
\section{The choice of automorphisms with incidence matrices of quadratic polynomials}

Assume the companion matrix related to a quadratic polynomial
$x^2 - ax \mp 1$ which is denoted by $A_{\pm}$:
\[
A_{\pm}=\pmatrix{ 0 & \pm 1  \cr 1 & a}
\]
is {\it hyperbolic}, that is,
the dominant eigenvalue $\lambda$ and the other
one $\lambda '$ hold $| \lambda |>1>|\lambda'|$,
then it is easily to check that there are following four cases.
\begin{proposition}
\label{prop:1-1} 
If a matrix $A_{-}=\pmatrix{ 0 & - 1  \cr 1 & a},~a \in
 {\mbox{\boldmath $Z$}}$ is hyperbolic, then there are two cases:\\
~~
\(
\begin{array}{ll}
{\rm (i)} & a \ge 3 \mbox{ and its eigenvalues } \lambda_1, \lambda_1' 
\mbox{ hold } 2 \le a-1 < \lambda_1 < a,
~0< \lambda_1'<1,\\
{\rm (ii)} & a \le -3 \mbox{ and its eigenvalues } \lambda_2, \lambda_2' 
\mbox{ hold } a < \lambda_2 < a+1 \le -2 ,~-1 < \lambda_2'<0.\\
\end{array}
\)\\
If a matrix $A_{+}=\pmatrix{ 0 &  1  \cr 1 & a},~a \in
 {\mbox{\boldmath $Z$}}$ is hyperbolic, then there are two cases:\\
~~\(
\begin{array}{ll}
{\rm (iii)} & a \ge 1 \mbox{ and its eigenvalues } \lambda_3, \lambda_3'
\mbox{ hold } 1 \le a < \lambda_3 < a+1,
~-1< \lambda_3'<0,\\
{\rm (iv)} & a \le -1 \mbox{ and its eigenvalues } \lambda_4, \lambda_4'
\mbox{ hold } a-1 < \lambda_4 < a \le -1,~
0 < \lambda_4'<1.\\
\end{array}
\)
\end{proposition}

As mention in Section 0, the aim of this paper is to find automorphisms
related to the matrices $A_{+}$ and $ A_{-}$ 
with which one can generate a stepped surface and
a Rauzy fractal, and discuss a dynamical system on the Rauzy fractal.
For this aim, the following automorphisms, which are
conjugate to some substitutions or some alternative substitutions, are
chosen for each case in Proposition \ref{prop:1-1}.
We call an endomorphism $\sigma$ on $F_2$ an
{\it ``alternative'' substitution} if {\it only} two letters $1^{-1}, 2^{-1}$ 
appear in $\sigma (i)$ for all $i \in {\cal A} $. If $\sigma$ is an alternative
substitution, then $\sigma ^{2}$ becomes a substitution.
That is the reason why such an endomorphism on $F_2$ is called an alternative
substitution.
We say an endomorphism $\sigma$ on $F_2$  is {\it conjugate} to 
an endomorphism $\tau$ if there exists an automorphism $\delta$ such that
$\sigma=\delta ^{-1} \circ \tau  \circ \delta $.

The case (i): A matrix is $A_{-}$ with $ a \ge 3$, 
and its eigenvalues hold $\lambda_1 > 1 ,~0<\lambda_1 ' < 1 $.
Set the automorphism $\sigma_1$ as 
\[
\sigma_1 :~
\left\{\
\begin{array}{l}
1 \to 2\\
2 \to 2^{a-2}1^{-1}22
\end{array}
\right.,
%%%%%%%%%
~~A_{{\sigma}_1}=A_{-}=
\left(
\begin{array}{cc}
0 & -1\\
1 & a
\end{array}
\right).
\]
The automorphism  $\sigma_1$
is conjugate to the invertible substitution $\tau_1$ with the automorphism
$\delta_1$ on $F_2$ such that
\[
\tau_1 :~
\left\{\
\begin{array}{l}
1 \to 2^{a-3}12\\
2 \to 2^{a-2}12
\end{array}
\right. ~,~
\delta_1 :~
\left\{\
\begin{array}{l}
1 \to 21^{-1}\\
2 \to 2
\end{array}
\right. ~~
\left(
\delta_1^{-1} :~
\left\{\
\begin{array}{l}
1 \to 1^{-1}2\\
2 \to 2
\end{array}
\right. ~~
\right).
\]

The case (ii): A matrix is $A_{-}$ with 
 $a \le -3$, and its eigenvalues hold 
$\lambda_2 < -1 ,~-1<\lambda_2 ' < 0 $.
Set the automorphism $\sigma_2$ as
\[
\sigma_2 :~
\left\{\
\begin{array}{l}
1 \to 2\\
2 \to 1^{-1} 2^{a}
\end{array}
\right. ,
%%%%%%%%%
~~A_{{\sigma}_2}=A_{-}=
\left(
\begin{array}{cc}
0 & -1\\
1 & a
\end{array}
\right).
\]
The automorphism  $\sigma_2$
is conjugate to the alternative substitution $\tau_2$ with the automorphism
$\delta_2$ on $F_2$ such that
\[
\tau_2 :~
\left\{\
\begin{array}{l}
1 \to 1^{-1}2^{a+2}\\
2 \to 1^{-1}2^{a+1}
\end{array}
\right. ~,~
\delta_2 :~
\left\{\
\begin{array}{l}
1 \to 2^{-1}1\\
2 \to 2
\end{array}
\right. ~~
\left(
\delta_2^{-1} :~
\left\{\
\begin{array}{l}
1 \to 21\\
2 \to 2
\end{array}
\right. ~~
\right).
\]

The case (iii): A matrix is $A_{+}$ with
$a \ge 1 $, and its eigenvalues hold
$\lambda_3 >1 ,~-1<\lambda_3 ' < 0 $.
Set the automorphism $\sigma_3$ as
\[
\sigma_3 :~
\left\{\
\begin{array}{l}
1 \to 2\\
2 \to 2^{a} 1
\end{array}
\right. ,
%%%%%%%%%
~~A_{{\sigma}_3}=A_{+}=
\left(
\begin{array}{cc}
0 & 1\\
1 & a
\end{array}
\right).
\]
In this case, the automorphism $\sigma_3$ is a substitution. 
Since the property in the case of substitutions is known as we saw in
Section 0, so we don't deal this case in this paper.

The case (iv): A matrix is $A_{+}$ with 
 $a\le -1$, and the eigenvalues hold 
$\lambda_4 < -1 ,~0<\lambda_4 ' < 1 $.
Set the automorphism $\sigma_4$ as
\[
\sigma_4 :~
\left\{\
\begin{array}{l}
1 \to 2\\
2 \to 1 2^{a}
\end{array}
\right. ,
%%%%%%%%%
~~A_{{\sigma}_4}=A_{+}=
\left(
\begin{array}{cc}
0 & 1\\
1 & a
\end{array}
\right).
\]
The automorphism  $\sigma_4$
is conjugate to the alternative substitution $\tau_4$ with the automorphism
$\delta_4$ on $F_2$ such that
\[
\tau_4 :~
\left\{\
\begin{array}{l}
1 \to 2^{-1}\\
2 \to 1^{-1}2^a
\end{array}
\right.~,~ 
\delta_4 :~
\left\{\
\begin{array}{l}
1 \to 1^{-1}\\
2 \to 2
\end{array}
\right. 
~~({\delta}_4^{-1}={\delta}_4).
\]
In this paper, we use the following typical examples for each case in figures:\\
\begin{center}
\(
\begin{array}{lll}
\sigma_1 :~
\left\{\
\begin{array}{l}
1 \to 2\\
2 \to 21^{-1}22
\end{array}
\right.~
&
\tau_1 :~
\left\{\
\begin{array}{l}
1 \to 12\\
2 \to 212
\end{array}
\right. ~
&
\delta_1 :~
\left\{\
\begin{array}{l}
1 \to 21^{-1}\\
2 \to 2
\end{array}
\right. 
\\
\sigma_2 :~
\left\{\
\begin{array}{l}
1 \to 2\\
2 \to 1^{-1} 2^{-1}2^{-1}2^{-1}
\end{array}
\right.~
&
\tau_2 :~
\left\{\
\begin{array}{l}
1 \to 1^{-1}2^{-1}\\
2 \to 1^{-1}2^{-1}2^{-1}
\end{array}
\right. ~
&
\delta_2 :~
\left\{\
\begin{array}{l}
1 \to 2^{-1}1\\
2 \to 2
\end{array}
\right. ~
\\
\sigma_4 :~
\left\{\
\begin{array}{l}
1 \to 2\\
2 \to 1 2^{-1}
\end{array}
\right.~
&
\tau_4 :~
\left\{\
\begin{array}{l}
1 \to 2^{-1}\\
2 \to 1^{-1}2^{-1}
\end{array}
\right.~
& 
\delta_4 :~
\left\{\
\begin{array}{l}
1 \to 1^{-1}\\
2 \to 2
\end{array}
\right. 
\end{array}
\)
\end{center}

%% file: autorank2-3.tex
\section{Stepped surfaces}

In this section, we construct the stepped surface of $P_{\sigma}$,
$\sigma=\sigma_1, \sigma_2,\sigma_4$ by using the fact that $\sigma$ is
conjugate to some substitution or some alternative substitution.
Here, $\tau$ is used for a substitution or an alternative substitution, and
$\sigma$ for an endomorphism on the free group $F_2$ of rank 2.
First let us consider the stepped surface of $P_\tau$,
$\tau=\tau_1,\tau_2,\tau_4$.
Notice that from the property of conjugate, the eigenvalues 
of $A_{\tau_t}$, $t=1,2,4$ are the same as the eigenvalues $\lambda_t,
\lambda_t '$ of $A_{\sigma_t}$.
The matrices  $A_{\tau_1},-A_{\tau_2},-A_{\tau_4}$ are
primitive, so each incidence matrix $A_{\tau_t}$ of $\tau_t $, $t=1,2,4$ 
has a positive column eigenvector ${\hako u}_{\tau_t}$
and a positive low eigenvector ${\hako v}_{\tau_t}$
corresponding to each eigenvalue
${\lambda}_1>1,~{\lambda}_2<-1,~{\lambda}_4<-1$ by Perron-Frobenius Theorem. 
When we consider arbitrary substitution or
alternative substitution $\tau$, 
assume that it satisfies the
hyperbolic, irreducible, unimodular 
conditions and $A_\tau$ or $-A_\tau$ is primitive hereafter.
For simplicity,
set the low eigenvector of $A_{\tau}$ as ${\hako v}_{\tau}=(1,\beta)$ with
some $\beta>0$ corresponding to the eigenvalue $\lambda_{\tau}$ with
$|\lambda_{\tau}|>1$. The stepped surfaces 
$ {\cal S}_{\tau}, {\cal S}_{\tau}'$ of 
the contractive eigenspace of $A_{\tau}$, which is given by
$P_{\tau}=\{ {\hako x}\in {\mbox{\boldmath $R$}}^2 \mid <{\hako
x},{}^t{\hako v}_{\tau}>=0\}$, 
are defined
analogously as in the case of substitutions as follows:
\begin{eqnarray*}
 {\cal S}_{\tau}&:=& \bigcup_{({\hako x}, i^*) \in S_{\tau}} ({\hako x}, i^*),\\
{\cal S}_{\tau}'&:=& \bigcup_{({\hako x}, i^*) \in S_{\tau}'} ({\hako x}, i^*),
\end{eqnarray*}
where
\begin{eqnarray}
\label{stepped1}
 S_{\tau} &:=& \left\{({\hako x}, i^*) \in  {\mbox{\boldmath $Z$}}^2 \times \{1^*,2^*\}  
\mid 
\left<{\hako x},{}^t{\hako v}_{\tau} \right> > 0,  
\left<{\hako x}-{\hako e}_i ,{}^t{\hako v}_{\tau}\right> \le
0\right\},\\
\label{stepped2}
 S_{\tau}'&:=& \left\{({\hako x}, i^*) \in  {\mbox{\boldmath $Z$}}^2 \times \{1^*,2^*\}  
\mid 
\left<{\hako x},{}^t{\hako v}_{\tau}\right> \ge 0,  
\left<{\hako x}-{\hako e}_i ,{}^t{\hako v}_{\tau}\right> < 0\right\}.
\end{eqnarray}
We mean by $({\hako x}, i^*)$
the positively oriented unit segment translated by ${\hako x}$
in  ${\mbox{\boldmath $Z$}}^2$, that is, 
\[
 ({\hako x},1^*):=\{{\hako x}+ t {\hako e}_2 \mid 0 \le t \le 1 \},~
 ({\hako x},2^*):=\{{\hako x}+ t {\hako e}_1 \mid 0 \le t \le 1 \}.
\] 
\begin{figure}[hbtp]
\setlength{\unitlength}{1mm}
\begin{center}
\begin{picture}(80,15)(0,0)
%\put(10,5){\line(0,1){10}}
%\put(0,15){\line(1,0){10}}
\put(0,5){\circle*{1}}
\put(0,5){\vector(0,1){10}}
\put(-5,0){$({\hako x}, 1^*)$}
\put(-5,5){${\hako x}$}
\put(20,5){\circle*{1}}
\put(20,15){\vector(0,-1){10}}
\put(15,0){$-({\hako x}, 1^*)$}
\put(15,5){${\hako x}$}
\put(40,5){\circle*{1}}
\put(40,5){\vector(1,0){10}}
\put(35,0){$({\hako x}, 2^*)$}
\put(35,5){${\hako x}$}
\put(60,5){\circle*{1}}
\put(70,5){\vector(-1,0){10}}
\put(55,0){$-({\hako x}, 2^*)$}
\put(55,5){${\hako x}$}
\end{picture}
\end{center}
\caption{The segments $({\hako x}, 1^*),({\hako x}, 2^*)$ with orientation}\label{fig:3}
\end{figure}

Notice that if $\beta$ is irrational, 
then
\[
 {\cal S} \setminus {\cal S}'=
\left\{
({\hako e}_1, 1^*) \cup ({\hako e}_2, 2^*)
\right\}
\setminus
\left\{
({\hako o}, 1^*) \cup ({\hako o}, 2^*)
\right\}.
\]
\begin{definition}
For an endomorphism $\sigma$ on $F_2$ given by
\[
 \sigma (i) = w_1^{(i)} w_2^{(i)} \cdots w_{l^{(i)}}^{(i)},~i \in {\cal A},
\]
define the $k$-prefix $ P_k^{(i)}$ and $k$-suffix  $S_k^{(i)} \in
F_2$  for $0 \le k \le l^{(i)}$  by 
\[
 P_k^{(i)} := w_1^{(i)} w_2^{(i)} \cdots w_{k-1}^{(i)},~
 S_k^{(i)} := w_{k+1}^{(i)} w_{k+2}^{(i)} \cdots w_{l^{(i)}}^{(i)}.
\]
Sometimes these notations are used for a substitution or an alternative
 substitution $\tau$ instead of $\sigma$.
The free ${\mbox{\boldmath $Z$}}$-module ${\cal G}^*$ is defined by
\[
{\cal G}^*:=\left\{ \sum_{k=1}^l n_k ({\hako x}_k, i_k^*) \mid 
n_k \in {\mbox{\boldmath $Z$}},~ 
{\hako x}_k \in {\mbox{\boldmath $Z$}}^2,~
i_k \in {\cal A} \mbox{ for any } k,~l < \infty \right\}, 
\]
whose element is identified with
a union of oriented unit segments with their multiplicity.
The tiling substitution $\sigma^*$ for a unimodular endomorphism $\sigma$ on $F_2$
such that $\det(A_{\sigma})=\pm 1$ is defined by 
\[
 \sigma^*({\hako x}, i^*):= \sum_{j \in {\cal A}}
\left\{
\sum_{w_k^{(j)} =i} \left(A_{\sigma}^{-1} ({\hako x}+ {\bf f}(S_k^{(j)})),j^*\right)
+
\sum_{w_k^{(j)} =i^{-1}} -\left(A_{\sigma}^{-1} ({\hako x}+ {\bf
f}(w_k^{(j)}S_k^{(j)})),j^*\right)
\right\}.
\]
\end{definition}
\begin{remark}
\label{remark:2-1}
In general, for a unimodular endomorphism $\sigma$ on the free group 
$F_d$ of rank $d$,  
a higher dimensional extension $E_k (\sigma)$ of $\sigma$ is defined
for $0 \le k \le d$, and $E_k^*(\sigma)$ is determined
as its dual map. The tiling substitution $\sigma^*$ is just
$E_1^*(\sigma)$ (cf. \cite{ei, s-a-i}).
\end{remark}
Define the subsets of ${\cal G}^*$ for a substitution or an alternative
substitution $\tau$ by
\begin{eqnarray}
\label{g1}
{{\cal G}_{\tau}^*} &:=& \left\{ \sum_{k=1}^l n_k ({\hako x}_k, i_k^*) \mid
\begin{array}{l}
n_k \in \{-1,1\},~ ({\hako x}_k, i_k^*) \in {S_{\tau}}, l<\infty\\
({\hako x}_k, i_k^*) \ne ({\hako x}_{k'}, i_{k'}^*) \mbox{ if } k \ne k'
\end{array}
\right\}, %\\
%\label{g2}
%{{\cal G}_{\tau}^*} '&:=& \left\{ \sum_{k=1}^l n_k ({\hako x}_k, i_k^*) \mid 
%\begin{array}{l}
%n_k \in \{-1,1\},~ ({\hako x}_k, i_k^*) \in {S_{\tau}}, l<\infty\\
%({\hako x}_k, i_k^*) \ne ({\hako x}_{k'}, i_{k'}^*) \mbox{ if } k \ne k'
%\end{array}
%\right\}, 
\end{eqnarray}
and ${{\cal G}_{\tau}^*} '$ is defined by replacing $ {S_{\tau}}$ with
$ S_{\tau}'$ in the formula (\ref{g1}).
For an element 
$\sum_{k=1}^l n_k({\hako x}_k, i_k^*) \in {{\cal G}_{\tau}^*}$,
the condition $n_k \in \{-1,1\}$ means that there is no overlap in it,
and we identify it with 
$\cup_{k=1}^l ({\hako x}_k, i_k^*) \subset {\cal S}_{\tau}$ geometrically.
The following two lemmas show that a tiling substitution $\tau^*$ 
is well-defined as a map on ${{\cal G}_{\tau}^*} $ 
(resp. a map from ${{\cal G}_{\tau}^*} $ to ${{\cal G}_{\tau}^*}'$)
for a substitution (resp. an alternative substitution) $\tau$.
\begin{lemma}
\label{lemma:2-2}
If $\tau$ is a substitution or an alternative substitution, then
$({\hako x}_1,i_1^*), ({\hako x}_2,i_2^*) \in S_{\tau}$,
$({\hako x}_1,i_1^*) \ne ({\hako x}_2,i_2^*)$  
implies $\tau^*({\hako x}_1,i_1^*) \cap
 \tau^*({\hako x}_2,i_2^*) = \emptyset$, where
 $\sum_{k_1=1}^{l_1}n_{k_1}({\hako
 x}_{k_1},i_{k_1}^*) \cap \sum_{k_2=1}^{l_2}
n_{k_2}({\hako x}_{k_2},i_{k_2}^*) \ne \emptyset$ means there exist 
$k_1 \in \{1,2,\cdots,l_1 \}$ and $k_2 \in \{1,2,\cdots,l_2 \}$ such that
$({\hako x}_{k_1},i_{k_1}^*) = ({\hako x}_{k_2},i_{k_2}^*)$.
\end{lemma}

Proof. 
In the case of substitutions, see \cite{a-i}.
We prove it for an alternative substitution. Suppose
$({\hako x}_1,i_1^*) , ({\hako x}_2,i_2^*) \in S_{\tau}$,
$({\hako x}_1,i_1^*) \ne ({\hako x}_2,i_2^*)$
 and $\tau^*({\hako x}_1,i_1^*) \cap
 \tau^*({\hako x}_2,i_2^*) \ne \emptyset$, then there exists
$j \in {\cal A} $, $k_1 \in \{1,2,\cdots, l^{(i_1)}\}$, $k_2 \in \{1,2,\cdots,
l^{(i_2)}\}$  
such that
$w_{k_1}^{(j)} =i_1^{-1}$, 
$w_{k_2}^{(j)} =i_2^{-1}$, 
and
\[
\left(A_{\tau}^{-1} ({\hako x_1}+ {\bf f}(w_{k_1}^{(j)}S_{k_1}^{(j)})),
j^*\right)
=
\left(A_{\tau}^{-1} ({\hako x_2}+ {\bf f}(w_{k_2}^{(j)} S_{k_2}^{(j)})),j^*\right).
\]
So 
\begin{eqnarray*}
A_{\tau}^{-1} ({\hako x_1}+ {\bf f}(w_{k_1}^{(j)} S_{k_1}^{(j)}))
&=&
A_{\tau}^{-1} ({\hako x_2}+ {\bf f}(w_{k_2}^{(j)} S_{k_2}^{(j)}))\\
{\hako x_1}+ {\bf f}(S_{k_1-1}^{(j)})
&=&
{\hako x_2}+ {\bf f}(S_{k_2-1}^{(j)})
\end{eqnarray*}
Suppose ${\hako x}_1 = {\hako x}_2$, 
then $k_1=k_2$ because $\tau$ is
an alternative substitution; and $i_1=i_2$. 
It contradicts to $({\hako x}_1,i_1^*) \ne ({\hako x}_2,i_2^*)$.
Therefore ${\hako x}_1 \ne {\hako x}_2$, and so $k_1 \ne k_2$.
We can suppose $k_1 < k_2$  without loss of generality, then
\begin{eqnarray*}
{\hako x_1} -{\hako e}_{i_1}
&=&
{\hako x}_2 + {\bf f}(S_{k_2-1}^{(j)}) - {\bf f}(S_{k_1-1}^{(j)})-{\hako
e}_{i_1}\\
&=&
{\hako x}_2 - {\bf f}(w_{k_1+1}^{(j)} \cdots  w_{k_2-1}^{(j)})
\end{eqnarray*}
and
\begin{eqnarray*}
<{\hako x}_1-{\hako e}_{i_1}, {}^t{\hako v}_{\tau}>
&=&
<{\hako x}_2 - {\bf f}(w_{k_1+1}^{(j)} \cdots  w_{k_2-1}^{(j)}) ,
{}^t{\hako v}_{\tau}>\\
&=&
<{\hako x}_2 ,{}^t{\hako v}_{\tau}> +
< - {\bf f}(w_{k_1+1}^{(j)} \cdots  w_{k_2-1}^{(j)}) ,{}^t{\hako
v}_{\tau}>>0
\end{eqnarray*}
It contradicts to $({\hako x}_1,i_1^*) \in {{\cal G}_{\tau}^*}$.
%In case where $k_1 > k_2$, we analogously show contradiction.
\hfill$\Box$ 
\\
%%%%%%%%%%
\begin{lemma}
\label{lemma:2-1}
If $\tau$ is a substitution, $({\hako x}, i^*) \in {{\cal G}_{\tau}^*}$
(resp. $({\hako x}, i^*) \in {{\cal G}_{\tau}^*}'$)
 implies $\tau^* ({\hako x}, i^*) \in {{\cal G}_{\tau}^*}$
 (resp. $\tau^* ({\hako x}, i^*) \in {{\cal G}_{\tau}^*}'$).
If $\tau$ is an alternative substitution, 
$({\hako x}, i^*) \in {{\cal G}_{\tau}^*}$
(resp. $({\hako x}, i^*) \in {{\cal G}_{\tau}^*}'$)
 implies $\tau^* ({\hako x}, i^*) \in {{\cal G}_{\tau}^*}'$
 (resp. $\tau^* ({\hako x}, i^*) \in {{\cal G}_{\tau}^*}$).
\end{lemma}

Proof. 
%For a substitution $\tau$, the result is known (cf. ???), so we prove it
%for a alternative substitution.
%In case where $\tau$ is a substitution,
%the tiling substitution is given by $E_1^*(\tau)({\hako x}, i^*)= \sum_{j=1}^2
%\sum_{w_k^{(j)} =i} \left(A_{\tau}^{-1} ({\hako x}+ {\bf
%f}(S_k^{(j)})),j^*\right)$.
%Suppose $({\hako x}, i^*) \in S_{\tau}$,
%then $({\hako x},{}^t{\hako v}_{\tau} ) > 0$, and  
%\( ({\hako x}-{\hako e}_i ,{}^t{\hako v}_{\tau}) \le 0 \) by
%the definition of $S_{\tau}$. 
%Because the Perron eigenvalue $\lambda_{\tau} >0$ and the low Perron
%eigenvector ${}^t{\hako v}_{\tau}>0$, if $w_k^{(j)} =i$ then
%\begin{eqnarray*}
%<A_{\tau}^{-1} ({\hako x}+ {\bf f}(S_k^{(j)})), {}^t{\hako v}_{\tau}> 
%&=&  
%<{\hako x}+ {\bf f}(S_k^{(j)}),{}^tA_{\tau}^{-1} {}^t{\hako v}_{\tau}>\\
%&=&  
%\frac{1}{\lambda_{\tau}}<{\hako x}+ {\bf f}(S_k^{(j)}),{}^t{\hako v}_{\tau}>>0,
%\end{eqnarray*}
%and
%\begin{eqnarray*}
%<A_{\tau}^{-1} ({\hako x}+ {\bf f}(S_k^{(j)}))-{\hako e}_j, {}^t{\hako
% v}_{\tau}>
%&=&  
%<{\hako x}+ {\bf f}(S_k^{(j)})-A_{\tau}{\hako e}_j,{}^tA_{\tau}^{-1} {}^t{\hako v}_{\tau}>\\
%&=&  
%\frac{1}{\lambda_{\tau}}\{<{\hako x} - {\hako e}_i,{}^t{\hako v}_{\tau}>
%+<- {\bf f}(P_k^{(j)}),{}^t{\hako v}_{\tau}>\}
%\le 0, 
%\end{eqnarray*}
%and we conclude $\left(A_{\tau}^{-1} ({\hako x}+ {\bf
%f}(S_k^{(j)})),j^*\right) \in S_{\tau}$ if $w_k^{(j)} =i$.
%
In the case of substitutions, see \cite{a-i}.
For an alternative substitution $\tau$,
the tiling substitution is given by $\tau^*({\hako x}, i^*)= \sum_{j \in
{\cal A}}
\sum_{w_k^{(j)} =i^{-1}} - \left(A_{\tau}^{-1} ({\hako x}+ {\bf
f}(S_k^{(j)}) -{\hako e}_i),j^* \right)$.
Suppose $({\hako x}, i^*) \in {{\cal G}_{\tau}^*}$.
Since the eigenvalue $\lambda_{\tau} <0$ 
and the eigenvector ${\hako v}_{\tau}>0$, if $w_k^{(j)} =i^{-1}$, then
\begin{eqnarray*}
<A_{\tau}^{-1} ({\hako x}+ {\bf f}(S_k^{(j)})-{\hako e}_i), {}^t{\hako v}_{\tau}> 
&=&  
<{\hako x}+ {\bf f}(S_k^{(j)})-{\hako e}_i,{}^tA_{\tau}^{-1} ~ {}^t{\hako v}_{\tau}>\\
&=&  
\frac{1}{\lambda_{\tau}}<{\hako x}+ {\bf f}(S_k^{(j)})-{\hako e}_i,{}^t{\hako
v}_{\tau}>\\
&=&  
\frac{1}{\lambda_{\tau}}
\{<{\hako x}-{\hako e}_i,{}^t{\hako v}_{\tau}>+
<{\bf f}(S_k^{(j)}),{}^t{\hako v}_{\tau}>\} \ge 0
\end{eqnarray*}
and
\begin{eqnarray*}
<A_{\tau}^{-1} ({\hako x}+ {\bf f}(S_k^{(j)})-{\hako e}_i)-{\hako
 e}_j, {}^t{\hako v}_{\tau}> 
&=&  
<{\hako x}+ {\bf f}(S_k^{(j)})-{\hako e}_i-A_{\tau}{\hako
e}_j,{}^tA_{\tau}^{-1}  ~ {}^t{\hako v}_{\tau}>\\
&=&  
\frac{1}{\lambda_{\tau}}<{\hako x}- {\bf f}(P_k^{(j)}),{}^t{\hako v}_{\tau}><0, 
\end{eqnarray*}
and we conclude $\left(A_{\tau}^{-1} ({\hako x}+ {\bf
f}(S_k^{(j)})-{\hako e}_i),j^*\right) \in S_{\tau}'$ and
$\tau^*({\hako x}, i^*) \in {{\cal G}_{\tau}^*}'$.
\hfill$\Box$\\ 

Recall the following lemma and proposition.
\begin{lemma}
\label{lemma:2-3}
(\cite{ei})
For unimodular endomorphisms $\sigma,~\sigma '$ on $F_2$,
the tiling substitution for their concatenation $\sigma \circ \sigma'$
is given by
$$(\sigma \circ \sigma')^* ={\sigma'}^* \circ \sigma^*.$$
\end{lemma}
\begin{proposition}
\label{prop:2-1}
(\cite{e-i})
If a substitution $\tau$ is invertible, then 
$\tau^{*~n} ({\cal U})$, $\tau^{*~n} ({\cal U}')$  and  
$\tau^{*~n} ({\hako e}_i,i^*)$, 
$\tau^{*~n} ({\hako o},i^*),~i \in {\cal A} $
are geometrically connected.
\end{proposition}
Since $({\hako e}_1, 1^*),({\hako e}_2, 2^*) \in S_{\tau}$ and
$({\hako o}, 1^*),({\hako o}, 2^*) \in S_{\tau}'$
for a substitution or an alternative substitution $\tau$,
so ${\cal U}:=({\hako e}_1, 1^*)+({\hako e}_2, 2^*) 
\in {{\cal G}_{\tau}^*},~
{\cal U}':=({\hako o}, 1^*)+({\hako o}, 2^*) \in {{\cal G}_{\tau}^*} '$.
Even if $\tau$ is an alternative substitutions,   
$\tau^2$ is a substitution and $(\tau^2)^*=(\tau^*)^2$.
And it is easy to check the substitutions $\tau_1, \tau_2^2,\tau_4^2$ are
invertible.
Thus we have the following proposition by Lemma \ref{lemma:2-2},
 Lemma \ref{lemma:2-1} and Proposition \ref{prop:2-1}.
\begin{proposition}
\label{prop:2-2}
In the case of (i),
\[
 \tau_1^{*~n} ({\cal U}) \in {\cal G}_{\tau_1}^* ,~
 \tau_1^{*~n} ({\cal U}')\in {{\cal G}_{\tau_1}^*} ',~ 
n \in {\mbox{\boldmath $N$}},
\]
and $ \tau_1^{*~n} ({\cal U}),~\tau_1^{*~n} ({\cal U}')$ are
connected.

In the case of (ii) and (iv),
\[
 \tau^{*~2n} ({\cal U})  \in {{\cal G}_{\tau}^*},~
 \tau^{*~2n} ({\cal U}') \in {{\cal G}_{\tau}^*}',~
n \in {\mbox{\boldmath $N$}},
\]
for $\tau = \tau_2, \tau_4$, and
$\tau^{*~2n} ({\cal U}'),~\tau^{*~2n} ({\cal U}) $ are connected.\\
\end{proposition}
%The proposition also says there is no overlap in $ E_1^{*~n}(\tau_1) ({\cal
% U}), E_1^{*~2n}(\tau) ({\cal U})$,
%$ E_1^{*~n}(\tau_1) ({\cal
% U}')$, $E_1^{*~2n}(\tau) ({\cal U}')$
% for $\tau=\tau_2,\tau_4$, 
%where we say there exists overlap in 
%$\sum_{k=1}^{l} n_k({\hako x}_k, i_k^*) \in {\cal G}^*$ 
%if $n_k \ge 2$ for some $k$.
\begin{remark}
\label{remark:2-2}
By using the idea of ${\cal C}$-covered property (cf. \cite{i-o, e-i-2}), we can show that
$\tau_1^{*~n} ({\cal U})$ 
(resp. $\tau_2^{*~2n} ({\cal U})$)  goes 
to the stepped surface ${\cal S}_{\tau_1}$ (resp. ${\cal S}_{\tau_2}$)
geometrically
 when $n$ goes to $\infty$.% It is true for ${\cal U}'$.
\end{remark}

The stepped surfaces ${\cal S}_{\tau},~{\cal S}_{\tau}'$ 
of the line $P_{\tau}$ for $\tau=\tau_1,\tau_2,\tau_4 $
are generated by using the tiling substitution with the seeds ${\cal U},~
{\cal U}'$.
%To discuss topological property of Rauzy fractals, let us observe the 
%topological property of $E_1^{*~n}(\tau_1) ({\cal U})$ and 
%$E_1^{*~2n}(\tau) ({\cal U})$ for $\tau=\tau_2, \tau_4$.
From now on we generate the stepped surface of the contractive
eigenspace $P_{\sigma}=
\{ {\hako x}\in {\mbox{\boldmath $R$}}^2 \mid 
<{\hako x} ,{}^t{\hako v}_{\sigma}>=0\}$, $\sigma=\sigma_1,\sigma_2,\sigma_4$ 
related to $A_{\sigma}$ in each case (i), (ii), (iv).
The matrices $A_{\sigma_i}$ and $-A_{\sigma_i}$ are not 
positive matrices, so we cannot apply Perron-Frobenius theorem
directly for them.
In fact, one of the eigenvalues of the incidence matrix $A_{\sigma_1}$ 
satisfies $\lambda_1>1$, 
and its corresponding low eigenvector ${\hako v}_{\sigma_1}$ 
given by $(1,\lambda_1)$ is positive, but
one of the eigenvalues of the incidence matrix $A_{\sigma_t}$ for
$t=2,4$ satisfies $\lambda_t<-1$,
and its corresponding low eigenvector ${\hako v}_{\sigma_t}$
given by $(-1,-\lambda_t)$ is not positive. 
So the sets $S_{\sigma_1},~S_{\sigma_1}'$ 
related to the stepped surface of 
the contractive eigenspace $P_{\sigma_1}$ are defined as
(\ref{stepped1}), (\ref{stepped2}), but
the sets $S_{\sigma},~S_{\sigma}',~\sigma=\sigma_2,\sigma_3$ 
are redefined as follows:
\begin{eqnarray*}
 S_{\sigma} &:=& \left\{({\hako x}, i^*) \in  {\mbox{\boldmath $Z$}}^2 \times \{1^*,2^*\}  
\left|
\begin{array}{ll}
<{\hako x},{}^t{\hako v}_{\sigma} > > 0,  
<{\hako x}+{\hako e}_1 ,{}^t{\hako v}_{\sigma}> \le 0
& \mbox{ if } i=1\\
<{\hako x}+{\hako e}_1,{}^t{\hako v}_{\sigma} > > 0,  
<{\hako x}+ {\hako e}_1-{\hako e}_2,{}^t{\hako v}_{\sigma}> \le 0
& \mbox{ if } i=2
\end{array}
\right.
\right\}
\\
 S_{\sigma}' &:=& \left\{({\hako x}, i^*) \in  {\mbox{\boldmath $Z$}}^2 \times \{1^*,2^*\}  
\left|
\begin{array}{ll}
<{\hako x},{}^t{\hako v}_{\sigma} > \ge 0,  
<{\hako x}+{\hako e}_1 ,{}^t{\hako v}_{\sigma}> < 0
& \mbox{ if } i=1\\
<{\hako x}+{\hako e}_1,{}^t{\hako v}_{\sigma} > \ge 0,  
<{\hako x}+{\hako e}_1-{\hako e}_2 ,{}^t{\hako v}_{\sigma}> < 0
& \mbox{ if } i=2
\end{array}
\right.
\right\}.
\end{eqnarray*}
The subset ${\cal G}_{\sigma}^*,~{{\cal G}_{\sigma}^*}',
~\sigma=\sigma_1,\sigma_2,\sigma_4$ of
${\cal G}^*$ are defined in the same way as (\ref{g1}) by $S_\sigma$, $S_\sigma'$.

Suppose an automorphism $\sigma$ is conjugate to 
$\tau$ as $\sigma=\delta ^{-1}
\circ \tau \circ \delta$ with some automorphism $\delta$ on $F_2$.
In general $A_{\sigma}= A_{\delta}^{-1} A_{\tau} A_{\delta}$, and
the contractive eigenspace $P_{\sigma}$ is given by
\begin{eqnarray*}
P_{\sigma}&=&\{A_{\delta}^{-1} {\hako x} \in  
{\mbox{\boldmath $R$}}^2  \mid {\hako x} \in  P_{\tau}\}\\
&=&A_{\delta}^{-1} P_{\tau}.
\end{eqnarray*}
By Lemma \ref{lemma:2-3}, the tiling substitution $\sigma^*$ of
$\sigma$ is
$$\sigma^*=\delta^* \circ \tau^* \circ
(\delta^{-1})^*,$$
and moreover,
\begin{eqnarray}
\label{decompE}
\sigma^{*~n}=\delta^* \circ \tau^{*~n} \circ
(\delta^{-1})^*.
\end{eqnarray}
%Notice $E_1^{*~n}(\tau)(g) \in {\cal G}_\tau$ if $g \in {\cal G}_\tau$.
%To generate the stepped surface of $P_{\sigma}$ by the tiling
%substitution $\sigma^*$, we need observe 
%$E_1^{*~n}(\delta)(g)$ for $g \in {\cal G}_\tau$.
The following replacing method will be introduced to understand
the relation between 
the stepped surfaces $S_{\tau_t}$ and $S_{\sigma_t},~t=1,2,4$ 
by using $\delta_t^*$.
%%%%%%%%%%%%%%%%%%%%
%the eigenvalues of $A_{\sigma}$ and $A_{\tau}$ are the same. So we denote
%the Perron eigenvalue of these matrix by $\lambda$. 
%%%%%%%%%%%%%%%%%%%%
\\
\\
{\bf Replacement Method}\\
\\
Choose low eigenvectors 
${\hako v}_{\tau_1}=(1,\frac{\lambda_1}{\lambda_1-1})$, 
${\hako v}_{\tau_2}=(1,\frac{\lambda_2}{\lambda_2+1})$,
${\hako v}_{\tau_4}=(1,-{\lambda_4})$
of $A_{{\tau}_t}$, $t=1,2,4$, then we have the following lemmas.
\begin{lemma}
\label{lemma:2-4}
~~
\begin{enumerate}
\item If $({\hako x},1^*) \in {\cal G}_{\tau_1}^*$, then $({\hako x},2^*) \in 
{\cal G}_{\tau_1}^*$.
\item If $({\hako x},1^*) \in {\cal G}_{\tau_2}^*$, then $({\hako x}-{\hako
      e}_1+{\hako e}_2,2^*) \in {\cal G}_{\tau_2}^*$.
\end{enumerate}
\end{lemma}

Proof.
We prove the first statement. The second one is proved by the same way.
Suppose $({\hako x},1^*) \in {\cal G}_{\tau_1}^*$, then 
$<{\hako x},{\hako v}_{\tau_1}> > 0$ and 
$<{\hako x}-{\hako e}_1,{\hako v}_{\tau_1}> \le 0$.
\begin{eqnarray*} 
<{\hako x}-{\hako e}_2,{\hako v}_{\tau_1}> &=&
<{\hako x}-{\hako e}_1,{\hako v}_{\tau_1}> + 
<{\hako e}_1-{\hako e}_2,{\hako v}_{\tau_1}> \\
&=& <{\hako x}-{\hako e}_1,{\hako v}_{\tau_1}> + 
1-\frac{\lambda_1}{\lambda_1-1} <0
\end{eqnarray*}
Therefore $({\hako x},2^*) \in {\cal G}_{\tau_1}^*$.
\hfill$\Box$ 
\\

%Now the replacement method is introduced for each case.
For a $2 \times 2$ matrix $A$ and $({\hako x}, i^*) \in 
{\cal G}^*$, 
$A ({\hako x}, i^*)
:=\{ A{\hako x}+ A {\hako y} \mid {\hako y} \in ({\hako x}, i^*)\}$.
We get $\delta_t^* ({\hako x},i^*),~t=1,2,4,~i \in {\cal A}$ by the following 
replacement method.

\begin{figure}[hbtp]
%%%%%% case i
The case (i):
\setlength{\unitlength}{1mm}
\begin{center}
%%%% (x,1)
\begin{picture}(140,20)(0,0)
\put(0,5){\circle*{1}}
\put(0,5){\vector(0,1){10}}
\put(-5,0){$({\hako x}, 1^*)$}
\put(-5,5){${\hako x}$}
\put(20,10){$A_{\delta_1}^{-1}$}
\put(20,5){$\longrightarrow$}
\put(40,5){\circle*{1}}
\put(40,5){\vector(0,1){10}}
\put(30,0){$A_{\delta_1}^{-1}{\hako x}$}
\put(50,10){replacement}
\put(60,5){$\longrightarrow$}
\put(80,5){\circle*{1}}
\put(80,15){\vector(0,-1){10}}
\put(70,0){$A_{\delta_1}^{-1}{\hako x}$}
\put(90,15){translation}
\put(90,10){by ${\hako e_1}-{\hako e_2}$}
\put(90,5){$\longrightarrow$}
\put(120,5){\circle*{1}}
\put(120,15){\vector(0,-1){10}}
\put(110,0){$A_{\delta_1}^{-1}{\hako x}+{\hako e_1}-{\hako e_2}$}
\end{picture}
%%%%% (x,2)
\begin{picture}(140,20)(0,0)
\put(0,5){\circle*{1}}
\put(0,5){\vector(1,0){10}}
\put(-5,0){$({\hako x}, 2^*)$}
\put(-5,5){${\hako x}$}
\put(20,5){$\longrightarrow$}
\put(45,5){\circle*{1}}
\put(45,5){\vector(-1,1){10}}
\put(35,0){$A_{\delta_1}^{-1}{\hako x}$}
\put(60,5){$\longrightarrow$}
\put(80,5){\circle*{1}}
\put(80,5){\vector(0,1){10}}
\put(70,15){\vector(1,0){10}}
\linethickness{.1pt}
\put(80,5){\line(-1,1){10}}
\put(70,0){$A_{\delta_1}^{-1}{\hako x}$}
\put(90,5){$\longrightarrow$}
\put(120,5){\circle*{1}}
\put(120,5){\vector(0,1){10}}
\put(110,15){\vector(1,0){10}}
\put(110,0){$A_{\delta_1}^{-1}{\hako x}+{\hako e_1}-{\hako e_2}$}
\end{picture}
%%%%% union
\begin{picture}(140,20)(0,0)
\put(0,5){\circle*{1}}
\put(0,5){\vector(0,1){10}}
\put(0,5){\vector(1,0){10}}
\put(-10,0){$({\hako x}, 1^*)+({\hako x}, 2^*)$}
\put(-5,5){${\hako x}$}
\put(20,5){$\longrightarrow$}
\put(45,5){\circle*{1}}
\put(45,5){\vector(0,1){10}}
\put(45,5){\vector(-1,1){10}}
\put(35,0){$A_{\delta_1}^{-1}{\hako x}$}
\put(60,5){$\longrightarrow$}
\put(80,5){\circle*{1}}
\put(80,5){\vector(0,1){10}}
\put(70,15){\vector(1,0){10}}
\put(82,15){\vector(0,-1){10}}
\linethickness{.1pt}
%\put(80,5){\line(-1,1){10}}
\put(70,0){$A_{\delta_1}^{-1}{\hako x}$}
\put(90,5){$\longrightarrow$}
\put(120,5){\circle*{1}}
%\put(120,5){\vector(0,1){10}}
\put(110,15){\vector(1,0){10}}
\put(110,0){$A_{\delta_1}^{-1}{\hako x}+{\hako e_1}-{\hako e_2}$}
\end{picture}
\end{center}
~\\
\epsfxsize=7cm
\epsfbox{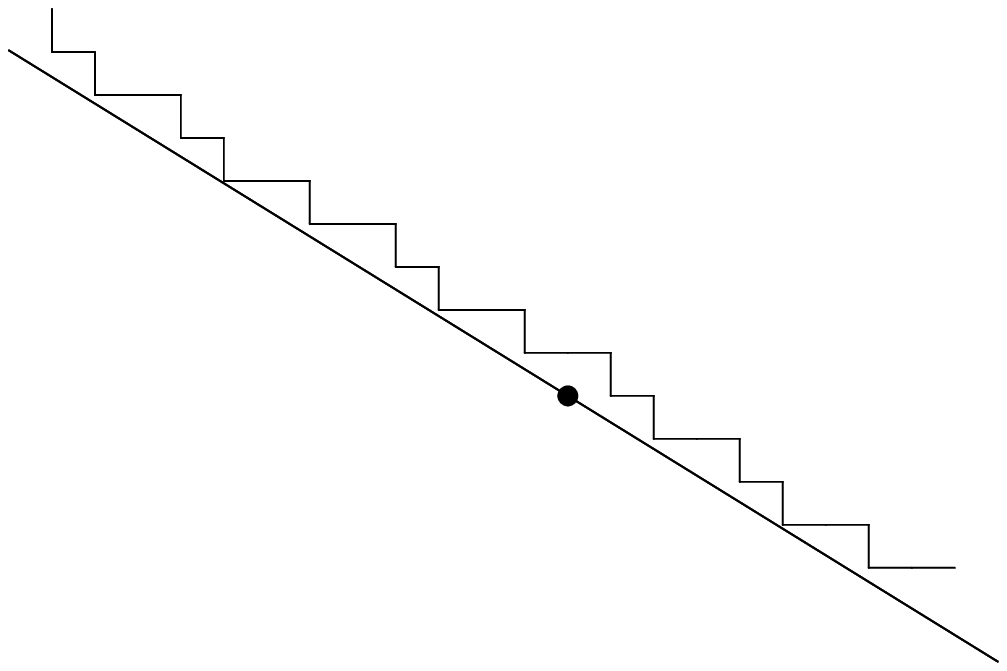}
~~~~~~~
\epsfxsize=7cm
\epsfbox{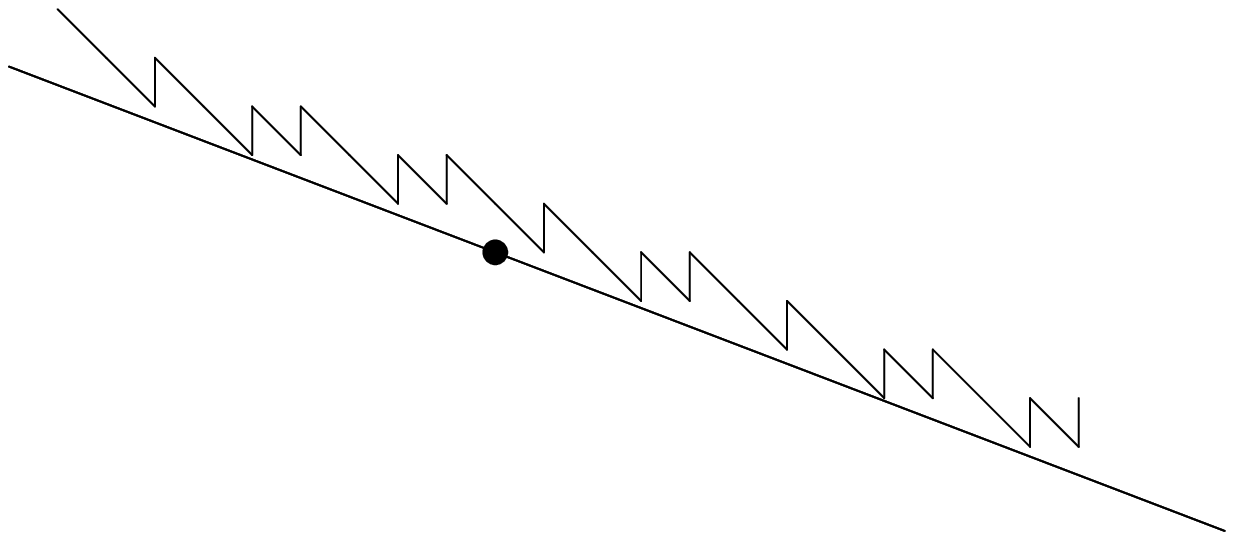}\\
\setlength{\unitlength}{1mm}
\begin{picture}(0,0)(0,0)
\put(0,40){$P_{\tau_1}$}
\put(70,30){$A_{\delta_1}^{-1}$}
\put(70,20){$\longrightarrow$}
\put(80,35){$P_{\sigma_1}$}
\end{picture}
The stepped surface ${\cal S}_{\tau_1}$
~~~~~~~~~~~~~~~~~~~~~~~~~~~
The picture after mapping by $A_{\delta_1}^{-1}$
\\
\vspace*{1cm}
\\
\epsfxsize=7cm
\epsfbox{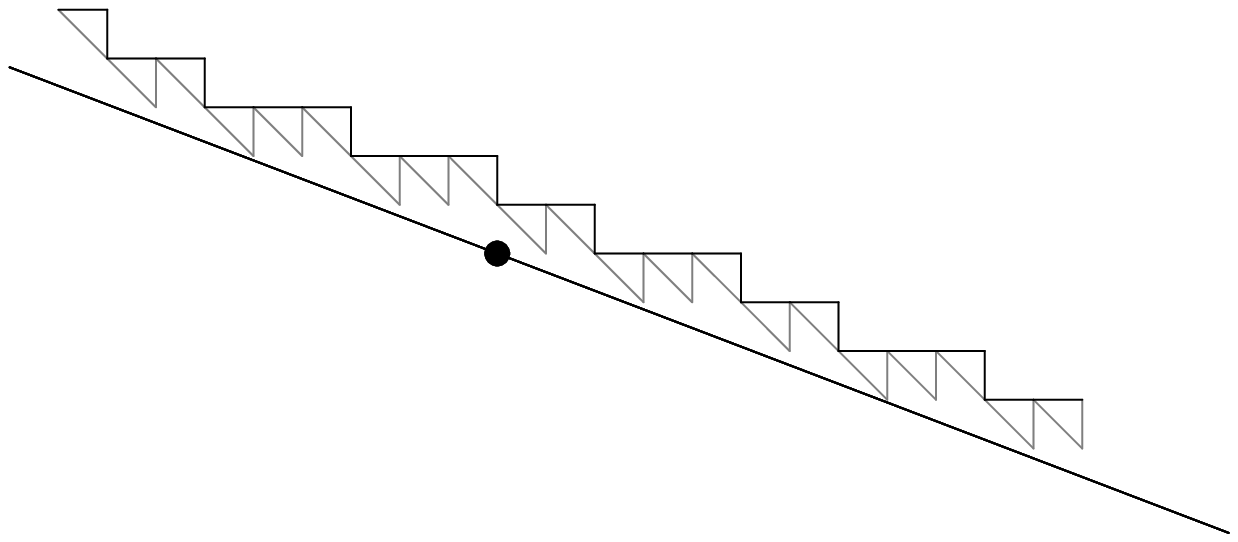}
~~~~~~~
\epsfxsize=7cm
\epsfbox{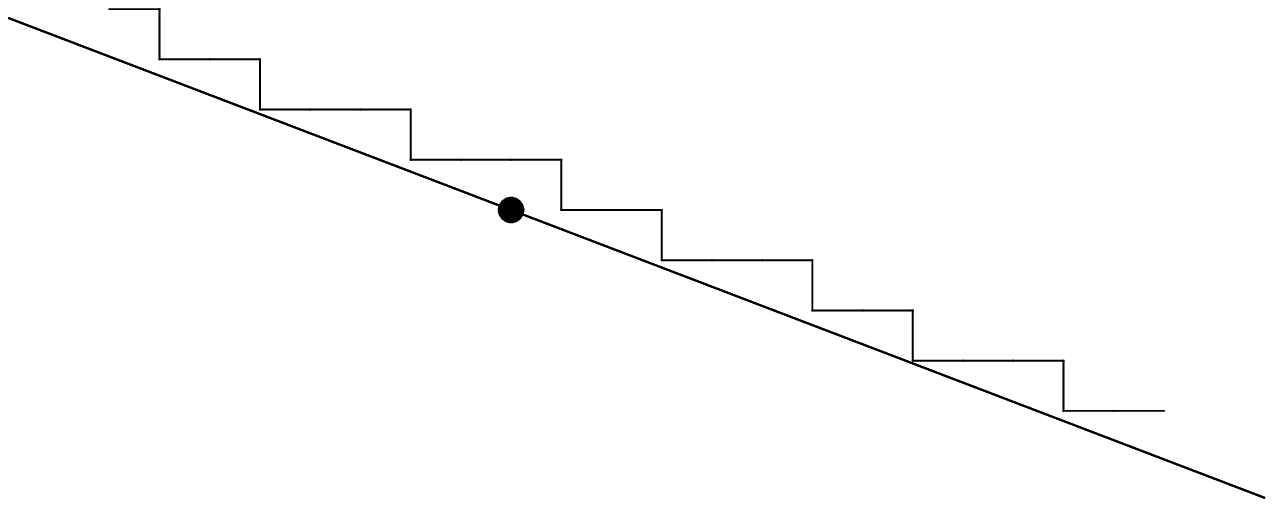}
\\
\setlength{\unitlength}{1mm}
\begin{picture}(0,0)(0,0)
\put(80,33){$P_{\sigma_1}$}
\end{picture}
The picture after replacement
~~~~~~~~~~~~~~~~~~~
The picture after translation by ${\hako e}_1 -{\hako e}_2$
\caption{Replacement and translation in the case (i)}\label{fig:4-1}
\end{figure}

In the case (i), for $({\hako x}, i^*) \in {\cal G}_{\tau_1}^*$
\[
 \delta_1^*({\hako x}, i^*)=\left\{
\begin{array}{ll}
-(A_{\delta_1}^{-1}{\hako x}+{\hako e}_1-{\hako e}_2,1^*) & i=1\\
(A_{\delta_1}^{-1}{\hako x}+{\hako e}_1-{\hako e}_2,1^*)+
(A_{\delta_1}^{-1}{\hako x},2^*)
 & i=2
\end{array}
\right. .
\]
Replace $A_{\delta_1}^{-1} ({\hako x}, 1^*)$ by 
$-(A_{\delta_1}^{-1} {\hako x}, 1^*)$, and translate it by 
${\hako e}_1-{\hako e}_2$, we get $\delta_1^* ({\hako x}, 1^*)$.
Replace $A_{\delta_1}^{-1} ({\hako x}, 2^*)$
by $(A_{\delta_1}^{-1} {\hako x}-{\hako e}_1+{\hako e}_2, 1^*)
+(A_{\delta_1}^{-1}{\hako x}, 2^*)$, 
and translate it by ${\hako e}_1-{\hako e}_2$, 
then we get  $\delta_1^* ({\hako x}, 2^*)$.
If $({\hako x}, 1^*) \in {\cal G}_{\tau_1}^*$, then
$({\hako x}, 1^*)+({\hako x}, 2^*)\in {\cal G}_{\tau_1}^*$ by Lemma
\ref{lemma:2-4}. 
Therefore the unit segment
$\delta_1^*({\hako x}, 1^*)=
-(A_{\delta_1}^{-1}{\hako x}+{\hako e}_1-{\hako e}_2,1^*)$
with negative orientation
is always cancelled by $\delta_1^*({\hako x}, 2^*)$  
(see figure \ref{fig:4-1}).\\

\begin{figure}[hbtp]
The case (ii):
\setlength{\unitlength}{1mm}
\begin{center}
%%%% (x,1)
\begin{picture}(140,20)(-20,0)
\put(0,5){\circle*{1}}
\put(0,5){\vector(0,1){10}}
\put(-5,0){$({\hako x}, 1^*)$}
\put(-5,5){${\hako x}$}
\put(20,10){$A_{\delta_2}^{-1}$}
\put(20,5){$\longrightarrow$}
\put(40,5){\circle*{1}}
\put(40,5){\vector(0,1){10}}
\put(30,0){$A_{\delta_2}^{-1}{\hako x}$}
\put(50,10){replacement}
\put(60,5){$\longrightarrow$}
\put(80,5){\circle*{1}}
\put(80,5){\vector(0,1){10}}
\put(70,0){$A_{\delta_2}^{-1}{\hako x}$}
\end{picture}
%%%%% (x,2)
\begin{picture}(140,20)(-20,0)
\put(0,5){\circle*{1}}
\put(0,5){\vector(1,0){10}}
\put(-5,0){$({\hako x}, 2^*)$}
\put(-5,5){${\hako x}$}
\put(20,5){$\longrightarrow$}
\put(35,5){\circle*{1}}
\put(35,5){\vector(1,1){10}}
\put(35,0){$A_{\delta_2}^{-1}{\hako x}$}
\put(60,5){$\longrightarrow$}
\put(70,5){\circle*{1}}
\put(80,15){\vector(0,-1){10}}
\put(70,5){\vector(1,0){10}}
\linethickness{.1pt}
\put(70,5){\line(1,1){10}}
\put(70,0){$A_{\delta_2}^{-1}{\hako x}$}
\put(90,5){$\longrightarrow$}
\put(110,5){\circle*{1}}
\put(120,15){\vector(0,-1){10}}
\put(110,5){\vector(1,0){10}}
\put(110,0){$A_{\delta_2}^{-1}{\hako x}$}
\end{picture}
%%%%% union
\begin{picture}(140,20)(-20,0)
\put(0,5){\circle*{1}}
\put(-10,15){\vector(1,0){10}}
\put(0,5){\vector(0,1){10}}
\put(-20,0){$({\hako x}, 1^*)+({\hako x}-{\hako e}_1+{\hako e}_2, 2^*)$}
\put(5,5){${\hako x}$}
\put(20,5){$\longrightarrow$}
\put(45,5){\circle*{1}}
\put(35,5){\vector(1,1){10}}
\put(45,5){\vector(0,1){10}}
\put(35,0){$A_{\delta_2}^{-1}{\hako x}$}
\put(60,5){$\longrightarrow$}
\put(80,5){\circle*{1}}
\put(80,15){\vector(0,-1){10}}
\put(82,5){\vector(0,1){10}}
\put(70,5){\vector(1,0){10}}
\put(70,0){$A_{\delta_2}^{-1}{\hako x}$}
\put(90,5){$\longrightarrow$}
\put(120,5){\circle*{1}}
\put(110,5){\vector(1,0){10}}
\put(110,0){$A_{\delta_2}^{-1}{\hako x}$}
\end{picture}
\end{center}
~\\
\epsfxsize=7cm
\epsfbox{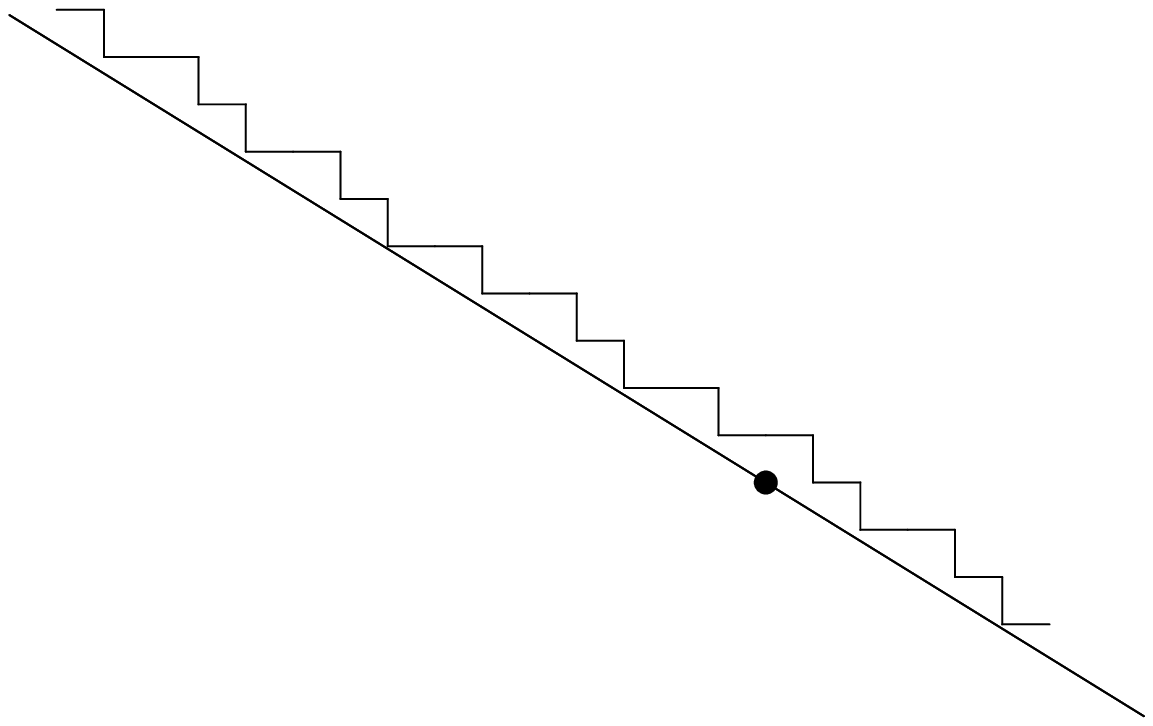}
~~~~~~~
\epsfxsize=7cm
\epsfbox{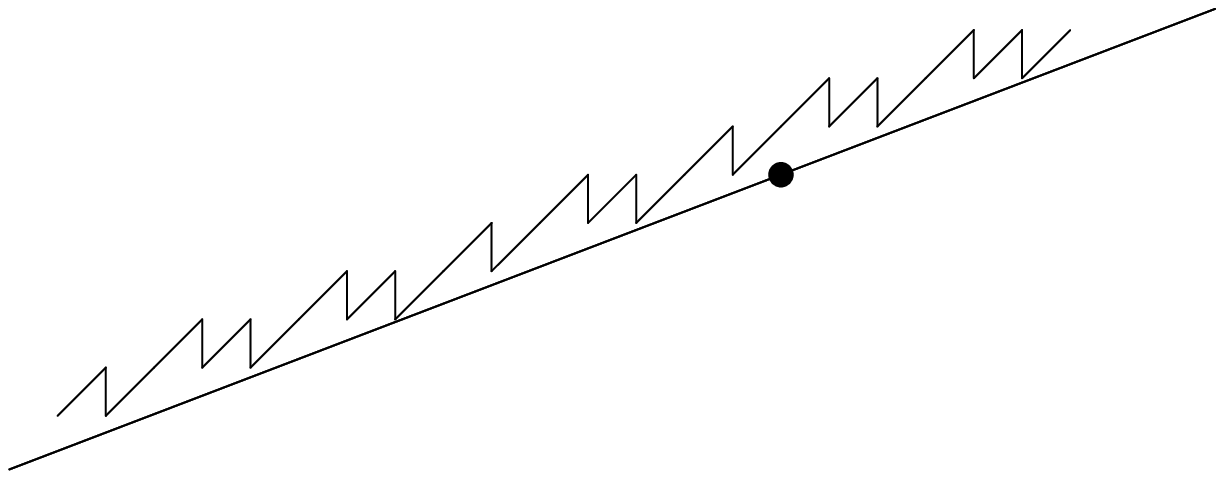}\\
\setlength{\unitlength}{1mm}
\begin{picture}(0,0)(0,0)
\put(0,40){$P_{\tau_2}$}
\put(70,30){$A_{\delta_2}^{-1}$}
\put(70,20){$\longrightarrow$}
\put(80,12){$P_{\sigma_2}$}
\end{picture}
The stepped surface ${\cal S}_{\tau_2}$
~~~~~~~~~~~~~~~~~~~~~~~~~~~
The picture after mapping by $A_{\delta_2}^{-1}$
\\
\vspace*{1cm}
\\
\epsfxsize=7cm
\epsfbox{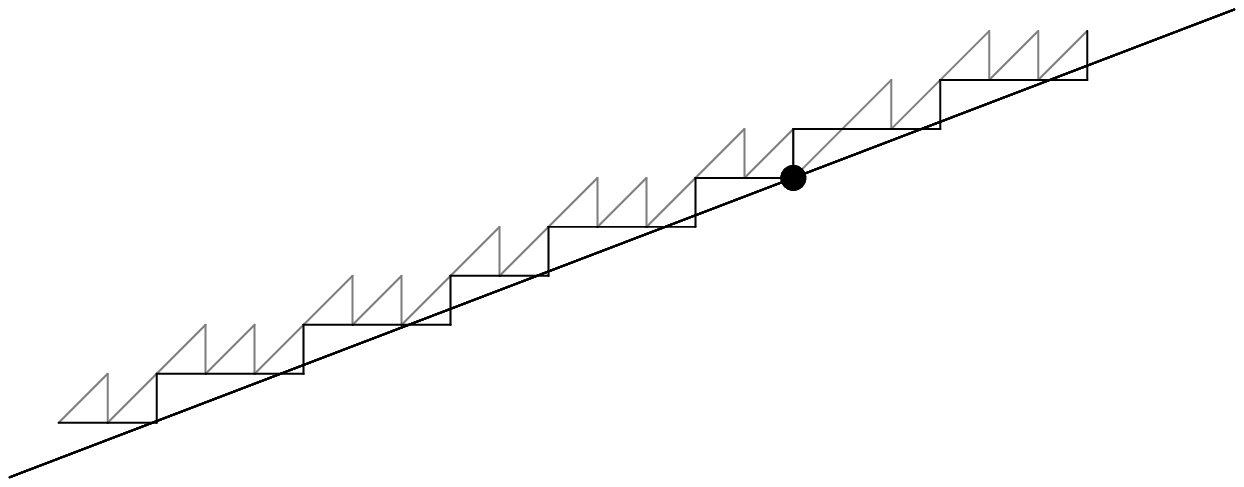}
~~~~~~~
\epsfxsize=7cm
\epsfbox{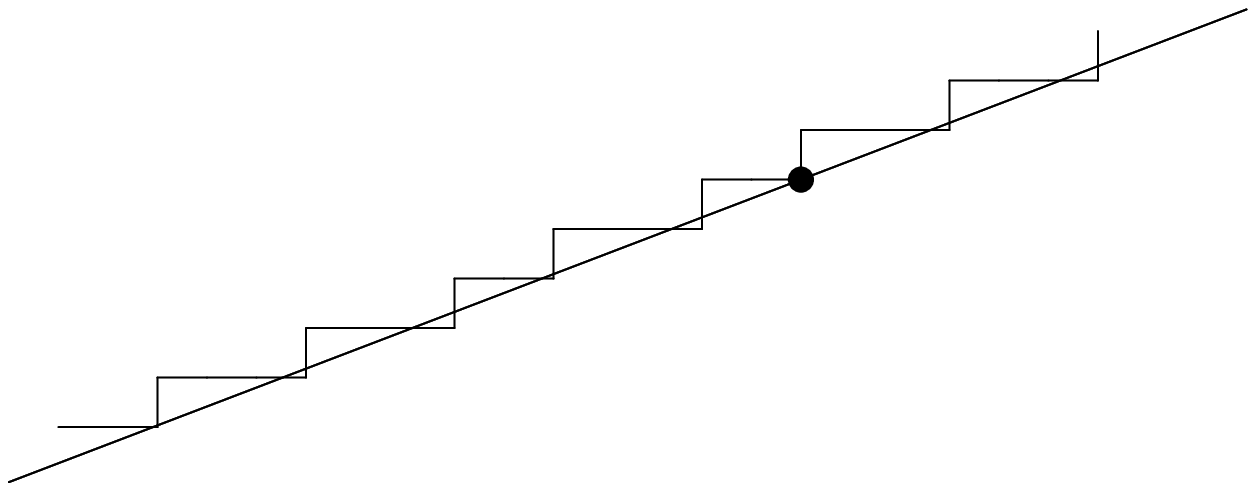}
\\
\setlength{\unitlength}{1mm}
\begin{picture}(0,0)(0,0)
\put(80,12){$P_{\sigma_2}$}
\end{picture}
The picture after replacement
~~~~~~~~~~~~~~~~~~~
%The picture after translation by ${\hako e}_1 -{\hako e}_2$
\caption{Replacement and translation in the case (ii)}\label{fig:4-2}
\end{figure}

In the case (ii), for $({\hako x}, i^*) \in {\cal G}_{\tau_2}^*$
\[
 \delta_2^*({\hako x}, i^*)=\left\{
\begin{array}{ll}
(A_{\delta_2}^{-1}{\hako x},1^*) & i=1\\
-(A_{\delta_2}^{-1}{\hako x}+{\hako e}_1,1^*)+
(A_{\delta_2}^{-1}{\hako x},2^*)
 & i=2
\end{array}
\right. .
\]
Replace $A_{\delta_2}^{-1} ({\hako x}, 1^*)$ by
$(A_{\delta_2}^{-1}{\hako x},1^*) $, 
then we get $\delta_2^* ({\hako x}, 1^*)$.
Replace $A_{\delta_2}^{-1} ({\hako x}, 2^*)$ by 
$-(A_{\delta_2}^{-1}{\hako x}+{\hako e}_1,1^*)+
(A_{\delta_2}^{-1}{\hako x},2^*)$,
then we get $\delta_2^* ({\hako x}, 2^*)$.
If $({\hako x}, 1^*) \in {\cal G}_{\tau_2}^*$, then
$({\hako x},1^*)+({\hako x}-{\hako e}_1+{\hako e}_2,2^*) 
\in {\cal G}_{\tau_2}^*$.
Therefore the unit segment 
$\delta_2^*({\hako x}, 1^*)=
(A_{\delta_2}^{-1}{\hako x},1^*)$ with positive orientation
is always cancelled by $\delta_2^*({\hako x}-{\hako e}_1+{\hako e}_2,
2^*)$
(see figure \ref{fig:4-2}).\\
%%%%%% case ii

\begin{figure}[hbtp]
The case (iv):
\setlength{\unitlength}{1mm}
\begin{center}
%%%% (x,1)
\begin{picture}(140,20)(-10,0)
\put(0,5){\circle*{1}}
\put(0,5){\vector(0,1){10}}
\put(-5,0){$({\hako x}, 1^*)$}
\put(-5,5){${\hako x}$}
\put(20,10){$A_{\delta_4}^{-1}$}
\put(20,5){$\longrightarrow$}
\put(40,5){\circle*{1}}
\put(40,5){\vector(0,1){10}}
\put(30,0){$A_{\delta_4}^{-1}{\hako x}$}
\put(50,10){replacement}
\put(60,5){$\longrightarrow$}
\put(80,5){\circle*{1}}
\put(80,15){\vector(0,-1){10}}
\put(70,0){$A_{\delta_4}^{-1}{\hako x}$}
\put(90,15){translation}
\put(90,10){by ${\hako e_1}$}
\put(90,5){$\longrightarrow$}
\put(120,5){\circle*{1}}
\put(120,15){\vector(0,-1){10}}
\put(110,0){$A_{\delta_4}^{-1}{\hako x}+{\hako e_1}$}
\end{picture}
%%%%% (x,2)
\begin{picture}(140,20)(-10,0)
\put(0,5){\circle*{1}}
\put(0,5){\vector(1,0){10}}
\put(-5,0){$({\hako x}, 2^*)$}
\put(-5,5){${\hako x}$}
\put(20,5){$\longrightarrow$}
\put(45,5){\circle*{1}}
\put(45,5){\vector(-1,0){10}}
\put(35,0){$A_{\delta_4}^{-1}{\hako x}$}
\put(60,5){$\longrightarrow$}
\put(80,5){\circle*{1}}
\put(70,5){\vector(1,0){10}}
\put(70,0){$A_{\delta_4}^{-1}{\hako x}$}
\put(90,5){$\longrightarrow$}
\put(120,5){\circle*{1}}
\put(120,5){\vector(1,0){10}}
\put(110,0){$A_{\delta_4}^{-1}{\hako x}$}
\end{picture}
\end{center}
~\\
\epsfxsize=7cm
\epsfbox{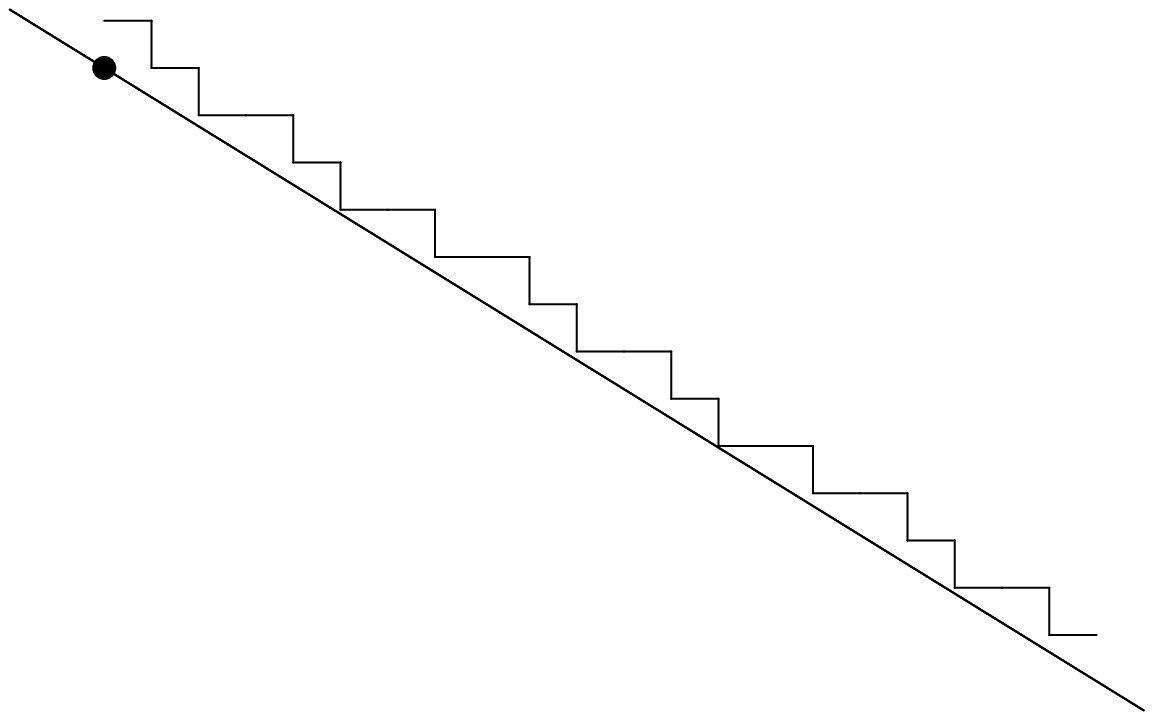}
~~~~~~~
\epsfxsize=7cm
\epsfbox{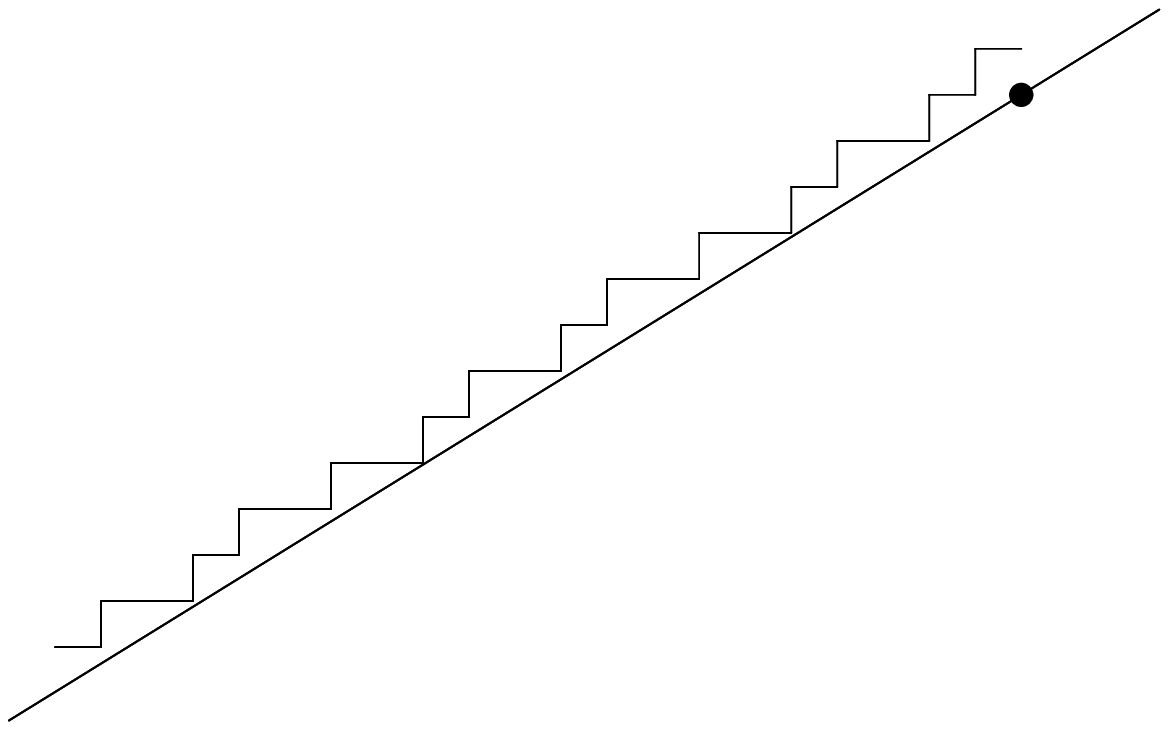}\\
\setlength{\unitlength}{1mm}
\begin{picture}(0,0)(0,0)
\put(0,40){$P_{\tau_4}$}
\put(70,30){$A_{\delta_4}^{-1}$}
\put(70,20){$\longrightarrow$}
\put(80,15){$P_{\sigma_4}$}
\end{picture}
The stepped surface ${\cal S}_{\tau_4}$
~~~~~~~~~~~~~~~~~~~~~~~~~~~
The picture after mapping by $A_{\delta_4}^{-1}$
\\
\vspace*{1cm}
\\
\epsfxsize=7cm
\epsfbox{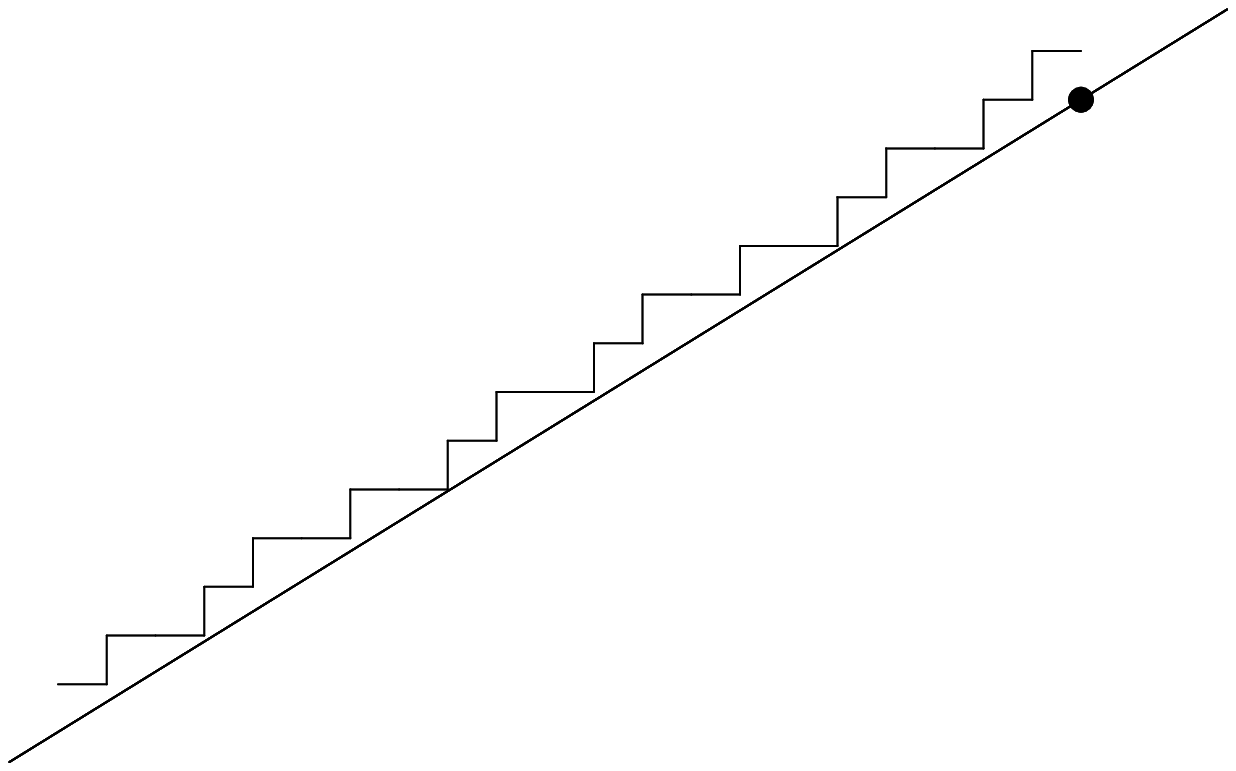}
~~~~~~~
\epsfxsize=7cm
\epsfbox{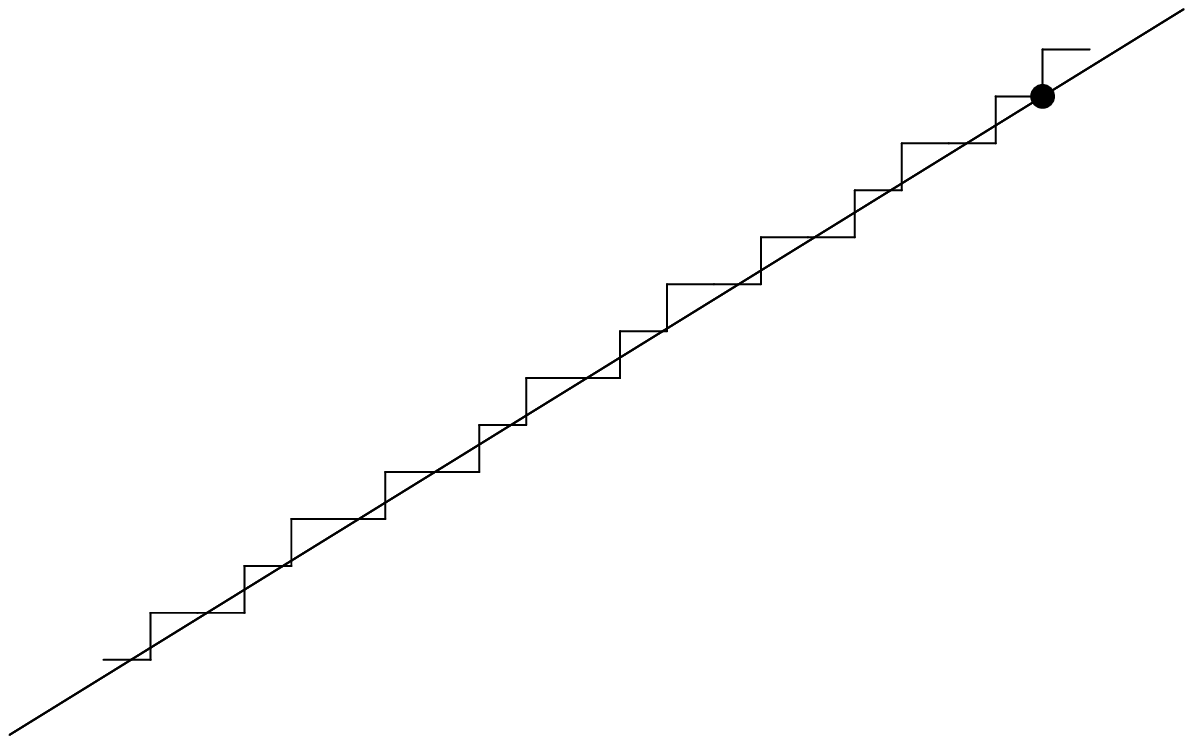}
\\
\setlength{\unitlength}{1mm}
\begin{picture}(0,0)(0,0)
\put(80,15){$P_{\sigma_4}$}
\end{picture}
The picture after replacement
~~~~~~~~~~~~~~~~~~~
The picture after translation by ${\hako e}_1 $
\caption{Replacement and translation in the case (iv)}\label{fig:4-4}
\end{figure}

In the case (iv), for $({\hako x}, i^*) \in {\cal G}_{\tau_4}^*$
\[
 \delta_4^*({\hako x}, i^*)=\left\{
\begin{array}{ll}
-(A_{\delta_4}^{-1}{\hako x}+{\hako e}_1,1^*) & i=1\\
(A_{\delta_4}^{-1}{\hako x},2^*) & i=2
\end{array}
\right. .
\]
Replace $A_{\delta_4}^{-1} ({\hako x}, 1^*)$ by
$-(A_{\delta_4}^{-1}{\hako x},1^*)$, and
translate it by ${\hako e}_1$, then
we get $\delta_4^* ({\hako x}, 1^*)$.
Replace $A_{\delta_4}^{-1} ({\hako x}, 2^*)$ by
$(A_{\delta_4}^{-1}{\hako x}-{\hako e}_1,2^*)$, and
translate it by ${\hako e}_1$, then
we get $\delta_4^* ({\hako x}, 2^*)$ (see figure \ref{fig:4-4}).
%%%%%% case iv

The following lemma shows the relation between the stepped surface
${\cal S}_{\tau_t}$ and ${\cal S}_{\sigma_t},~t=1,2,4$ by using $\delta_t^*$.
For ${\hako y} \in {\mbox{\boldmath $Z$}}^2$ and
$\sigma=\sigma_1,\sigma_2,\sigma_3$,
\[
{\cal G}_{\sigma}^*+{\hako y}:=
 \left\{\sum_{k=1}^l n_k({\hako x}_k + {\hako y}, i_k^*) \mid
\sum_{k=1}^l n_k({\hako x}_k, i_k^*) \in  {\cal G}_{\sigma}^*
\right\}.
\]
\begin{lemma}
\label{lemma:2-5}
~~
\begin{enumerate}
\item 
\(
\begin{array}{l}
\delta_1^*(({\hako x},1^*)+({\hako x},2^*)) \in {\cal G}_{\sigma_1}^* 
~{\mbox if}
({\hako x},1^*) \in {\cal G}_{\tau_1}^*\\
\delta_1^*({\hako x},2^*) \in {\cal G}_{\sigma_1}^* 
~{\mbox if}
({\hako x},2^*) \in {\cal G}_{\tau_1}^* {\mbox ~and~} ({\hako x},1^*) 
\notin {\cal G}_{\tau_1}^*
\end{array}
\)
\item 
\(
\begin{array}{l}
\delta_2^*(({\hako x},1^*)+({\hako x}-{\hako e}_1+{\hako e}_2,2^*)) 
\in {\cal G}_{\sigma_2}^*+{\hako e}_1 
~{\mbox if~}
({\hako x},1^*) \in {\cal G}_{\tau_2}^*\\
\delta_2^*({\hako x},2^*)
\in {\cal G}_{\sigma_2}^*+{\hako e}_1 
~{\mbox if~}
 ({\hako x},2^*) \in {\cal G}_{\tau_2}^* {\mbox ~and~}
({\hako x}+{\hako e}_1-{\hako e}_2,1^*) \notin {\cal G}_{\tau_2}^* 
\end{array}
\)
\item 
$\delta_4^*({\hako x},i^*) \in {\cal G}_{\sigma_4}^* +{\hako e}_1 $ 
if $({\hako x},i^*) \in {\cal G}_{\tau_4}^*$
\end{enumerate}
%By replacing ${\cal G}_{\sigma_i}$ and ${\cal G}_{\tau_i} $ by
%${\cal G}_{\sigma_i} '$ and ${\cal G}_{\tau_i}  '$, we have analogous result.
\end{lemma}

Proof.
In the case where $({\hako x},1^*) \in {\cal G}_{\tau_1}^*$,
$({\hako x},1^*)+({\hako x},2^*)  \in {\cal G}_{\tau_1}^*$ and 
$\delta_1^*(({\hako x},1^*)+({\hako x},2^*) ) =
(A_{\delta_1}^{-1}{\hako x},2^*)$.
From ${\hako v}_{\sigma_1} A_{\delta_1}^{-1} =(\lambda_1 -1)
{\hako v}_{\tau_1}$,
\begin{eqnarray*} 
<A_{\delta_1}^{-1}{\hako x},{}^t{\hako v}_{\sigma_1}>&=&
<{\hako x},{}^tA_{\delta_1}^{-1}~{}^t{\hako v}_{\sigma_1}>\\
&=& (\lambda_1 -1) <{\hako x},{}^t{\hako v}_{\tau_1}> >0,
\end{eqnarray*}
and
\begin{eqnarray*} 
<A_{\delta_1}^{-1}{\hako x}-{\hako e}_2,{}^t{\hako v}_{\sigma_1}>&=&
<{\hako x}-{\hako e}_2,{}^tA_{\delta_1}^{-1}~{}^t{\hako v}_{\sigma_1}>\\
&=& (\lambda_1 -1) <{\hako x}-{\hako e}_2,{}^t{\hako v}_{\tau_1}> \le 0.
\end{eqnarray*}
Therefore
$(A_{\delta_1}^{-1}{\hako x},2^*) \in  {\cal G}_{\sigma_1}^*$. 

In the case where 
$({\hako x},2^*)  \in {\cal G}_{\tau_1}^*$ and 
$({\hako x},1^*)  \notin {\cal G}_{\tau_1}^*$, 
%we will prove 
%$E_1^*(\delta_1)({\hako x},2^*) =
%(A_{\delta_1}^{-1}{\hako x}+{\hako e}_1-{\hako e}_2,1^*)+
%(A_{\delta_1}^{-1}{\hako x},2^*) \in  {\cal G}_{\sigma_1}$.
noticing
$<{\hako x}-{\hako e}_1,{}^t{\hako v}_{\tau_1} > > 0$
by $({\hako x},1^*)  \notin {\cal G}_{\tau_1}^*$, we have
\begin{eqnarray*} 
<A_{\delta_1}^{-1}{\hako x}+{\hako e}_1-{\hako e}_2,
{}^t{\hako v}_{\sigma_1}>
&=&
<{\hako x}-{\hako e}_1,{}^tA_{\delta_1}^{-1}~{}^t{\hako v}_{\sigma_1}>\\
&=& (\lambda_1 -1) <{\hako x}-{\hako e}_1,{}^t{\hako v}_{\tau_1}> > 0,\\
<A_{\delta_1}^{-1}{\hako x}-{\hako e}_2,
{}^t{\hako v}_{\sigma_1}>
&=&
 (\lambda_1 -1) <{\hako x}-{\hako e}_2,{}^t{\hako v}_{\tau_1}> \le 0,\\
<A_{\delta_1}^{-1}{\hako x},
{}^t{\hako v}_{\sigma_1}>
&=&
 (\lambda_1 -1) <{\hako x},{}^t{\hako v}_{\tau_1}> > 0.\\
\end{eqnarray*}
Therefore
$(A_{\delta_1}^{-1}{\hako x}+{\hako e}_1-{\hako e}_2,1^*)+
(A_{\delta_1}^{-1}{\hako x},2^*) \in  {\cal G}_{\sigma_1}^*$. 
The first statement is proved, and the others can
 be proved analogously.
\hfill$\Box$ \\

By the replacement method, we have the following lemma.
\begin{lemma}
\label{lemma:2-6}
If $\gamma \in {\cal G}_{\tau_t}^*$, $t=1,2,4$ is connected, 
% or $\gamma \in {\cal G}_{\tau_t}'$
then $\delta_t^* (\gamma)$ is also connected.
\end{lemma}

To generate the stepped surface of $P_{\sigma_t},~t=1,2,4$, determine an initial element
${\widetilde{\cal U}}$ and $\widetilde{\cal U}'$ 
for $\sigma_t^*$ as follows:
\[
\widetilde{\cal U} := \delta_t^* ({\cal U} ),~
\widetilde{\cal U}' := \delta_t^* ({\cal U} ').
\]
%(See figure???)
%From Proposition 2.1, Lemma 5, Lemma 6, we obtain the theorem.
\begin{theorem}
For any positive integer $n$, 
\begin{enumerate}
 \item $\sigma_1^{*~n} (\widetilde{\cal U}) \in {\cal G}_{\sigma_1}^*$,
%$\sigma_1^{*~n} (\widetilde{\cal U}') \in {\cal G}_{\sigma_1} ',$
\item $\sigma_t^{*~2n} (\widetilde{\cal U}) 
\in {\cal G}_{\sigma_t}^*+{\hako e}_1$,
%$E_1^{*~2n}(\sigma_2) (\widetilde{\cal U}')
%\in {\cal G}_{\sigma_t}' +{\hako e}_1$
$t=2,4$.
\end{enumerate}
Moreover, $\sigma_t^{*~n} (\widetilde{\cal U}),~t=1,2,4$ 
%E_1^{*~n}(\sigma_t) (\widetilde{\cal U}')$,$
are connected.
\end{theorem}

Proof.
From the equality (\ref{decompE}), 
\begin{eqnarray*}
\sigma_t^{*~n}(\widetilde{\cal U}) &=&
%E_1^{*}(\delta_i) \circ
%E_1^{*~n}(\tau_i) \circ  E_1^{*}(\delta_i^{-1})(\widetilde{\cal U})
%&=& 
\delta_t^* \circ
\tau_t^{*~n} ({\cal U}),~t=1,2,4.
\end{eqnarray*}
By Proposition \ref{prop:2-2}, Lemma \ref{lemma:2-5} and Lemma
\ref{lemma:2-6},
$\sigma_1^{*~n}(\widetilde{\cal U})$ 
(resp. $\sigma_t^{*~2n}(\widetilde{\cal U})$, $t=2,4$) is included
in ${\cal G}_{\sigma_1}^*$ (resp. ${\cal G}_{\sigma_t}^*+{\hako e_1}$), and
connected. The other cases can be proved analogously.
\hfill$\Box$ \\

\begin{figure}[hbtp]
\setlength{\unitlength}{1mm}
\begin{center}
\epsfxsize=2cm
\epsfbox{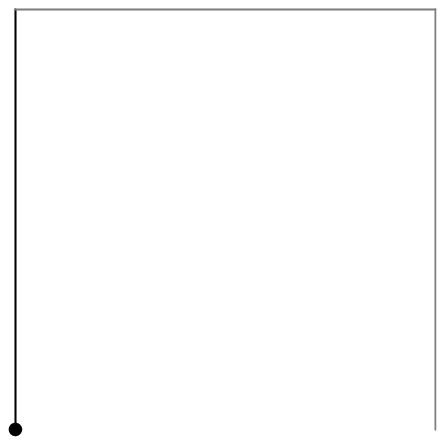}
~~~~
\epsfxsize=8cm
\epsfbox{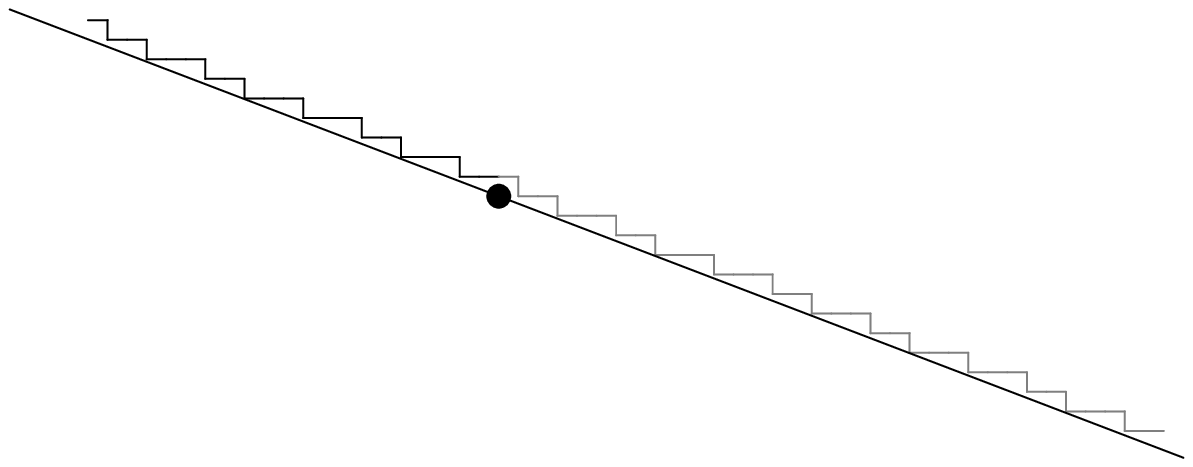}
\end{center}
\caption{The seed $\widetilde{\cal U}$ and $\sigma_1^{*~n}(\widetilde{\cal
 U})$ in the case (i)}\label{fig:5}
\end{figure}

By Remark \ref{remark:2-2}, $\sigma_1^{*~n}(\widetilde{\cal U})$
(resp. $\sigma_2^{*~2n}(\widetilde{\cal U})$) goes to the stepped
surface ${\cal S}_{\sigma_1}$ (resp. ${\cal S}_{\sigma_2}$) when $n$
goes to infinity (see Figure \ref{fig:5} and the last pictures of Figure
\ref{fig:4-1}, \ref{fig:4-2}, \ref{fig:4-4}).

%% file: autorank2-4.tex
\section{Rauzy fractals and domain exchange transformations}

In this section, we construct Rauzy fractals induced from automorphisms
$\sigma_t,~t=1,2,4$, and consider domain exchange transformations.

Define the projection $\pi_{\tau_t}$ (resp. $\pi_{\sigma_t}$),  
$t=1,2,4$ from ${\mbox{\boldmath $R$}}^2$ to
the contractive eigenspace $P_{\tau_t}$ (resp. $P_{\sigma_t}$)
along a column eigenvector ${\hako u}_{\tau_t}$ 
(resp. ${\hako u}_{\sigma_t}$) of $A_{\tau_t}$ (resp. $A_{\sigma_t}$)
corresponding to the eigenvalue $\lambda_t$.
First we define Rauzy fractals related to 
the substitution $\tau_1$ and  the alternative substitutions $\tau_2$, 
$\tau_4 $ as follows:
\begin{eqnarray*}
X_{\tau}&:=& \lim_{n \to \infty} A_{\tau}^n \pi_{\tau}
\tau^{*~n}({\cal U}),\\
&:=& \lim_{n \to \infty} A_{\tau}^n \pi_{\tau}
 \tau^{*~n}({\cal U}'),\\
X_{\tau}^{(i)}&:=& \lim_{n \to \infty} A_{\tau}^n \pi_{\tau}
 \tau^{*~n}({\hako e}_i, i^*),\\
{X'}_{\tau}^{(i)}  &:=& \lim_{n \to \infty} A_{\tau}^n \pi_{\tau}
 \tau^{*~n}({\hako o}, i^*),
\end{eqnarray*}
$\tau=\tau_1,\tau_2,\tau_4$.
It is proved that the limit sets exist in the sense of Hausdorff metric
by the same way in the case where $\tau$ is a substitution (cf. \cite{a-i}).
\begin{remark}
\label{remark:3-1}
Notice that one can replace $n$ with $2n$ in the formulas of definitions 
of Rauzy fractals above.
Thus for alternative substitutions $\tau=\tau_2, \tau_4$, $i \in {\cal A}$,
\[
 X_{\tau}^{(i)}=  X_{\tau^2}^{(i)},~
{X'}_{\tau}^{(i)} ={X'}_{\tau^2}^{(i)}.
\]
Therefore we can also apply Theorem \ref{theorem:0-2} 
and show that $X_{\tau}^{(i)},{X'}_{\tau}^{(i)}$ are intervals.
\end{remark}

\begin{proposition}
\label{prop:3-1}
A substitution or an alternative substitution $\tau$ is written as 
$\tau(i)=w_1^{(i)}w_2^{(i)}\cdots w_{l^{(i)}}^{(i)}$ and 
$\tau^2(i)=w_1^{(2,i)}w_2^{(2,i)}\cdots w_{l^{(2,i)}}^{(2,i)},~i \in {\cal A}$.
We denote by $P_k^{(i)}$ and $S_k^{(i)}$ (resp. $P_k^{(2,i)}$ and $S_k^{(2,i)}$)
the $k$-prefix and the $k$-suffix of $\tau(i)$  
(resp. $\tau^2(i)$).

In the case of (i),
\begin{eqnarray*}
A_{\tau_1}^{-1}X_{\tau_1}^{(i)}  &= & \cup_{j \in {\cal A}} \cup_{w_k^{(j)}=i} 
(-A_{\tau_1}^{-1} \pi_{\tau_1}{\bf f}(P_k^{(j)}) + X_{\tau_1}^{(j)} ),\\
A_{\tau_1}^{-1}{X'}_{\tau_1}^{(i)}  &= &\cup_{j \in {\cal A}} \cup_{w_k^{(j)}=i} 
(A_{\tau_1}^{-1}  \pi_{\tau_1}{\bf f}(S_k^{(j)}) + {X'}_{\tau_1}^{(j)} ).
\end{eqnarray*}

In the case of (ii) and (iv), for $\tau=\tau_2, \tau_4$,
\begin{eqnarray*}
A_{\tau^2}^{-1}X_{\tau}^{(i)}  &= & \cup_{j \in {\cal A}} \cup_{w_k^{(2,j)}=i} 
(-A_{\tau^2}^{-1}  \pi_{\tau}{\bf f}(P_k^{(2,j)}) + X_{\tau}^{(j)} ),\\
A_{\tau^2}^{-1}{X'}_{\tau}^{(i)}  &= &\cup_{j \in {\cal A}} \cup_{w_k^{(2,j)}=i} 
(A_{\tau^2}^{-1}  \pi_{\tau}{\bf f}(S_k^{(2,j)}) + {X'}_{\tau}^{(j)} ),
\end{eqnarray*}
and moreover, 
\begin{eqnarray*}
A_{\tau}^{-1}X_{\tau}^{(i)}  
%&= & \cup_{j=1}^2 \cup_{w_k^{(j)}=i^{-1}} 
%(A_{\tau}^{-1} {\bf f}(S_k^{(j)}) + {X'}_{\tau}^{(j)} ),\\
&= & \cup_{j \in {\cal A}} \cup_{w_k^{(j)}=i^{-1}} 
(-A_{\tau}^{-1}  \pi_{\tau}{\bf f}(P_k^{(j)}w_k^{(j)})+ {X}_{\tau}^{(j)} ),\\
A_{\tau}^{-1}{X'}_{\tau}^{(i)}  
%&= &\cup_{j=1}^2 \cup_{w_k^{(j)}=i^{-1}} 
%(-A_{\tau}^{-1} {\bf f}(P_k^{(j)}) + {X}_{\tau}^{(j)} )\\
 &= &\cup_{j \in {\cal A}} \cup_{w_k^{(j)}=i^{-1}} 
(A_{\tau}^{-1}  \pi_{\tau}{\bf f}(w_k^{(j)}S_k^{(j)})+ {X'}_{\tau}^{(j)} ).
\end{eqnarray*}
These unions are pairwise disjoint in the sense of Lebesgue measure. 
\end{proposition}

Proof.
For the substitutions $\tau_1, \tau^2_2,\tau^2_4$, these set equations
are known (see Proposition \ref{prop:0-3}). So we will show the last equations for alternative substitutions
$\tau_2,\tau_4$.
\begin{eqnarray*}
A_{\tau}^{-1}X_{\tau}^{(i)}  
&= & A_{\tau}^{-1} \lim_{n \to \infty} A_{\tau}^{n+1} \pi_{\tau}
 \tau^{*~(n+1)}({\hako e}_i, i^*)\\
&= & \lim_{n \to \infty} A_{\tau}^{n} \pi_{\tau}
 \tau^{*~n}
\left(
\sum_{j \in {\cal A}} \sum_{w_k^{(j)}=i^{-1}}
-(A_{\tau}^{-1}{\bf f}(S_k^{(j)}),j^*)
\right)\\
&= & \lim_{n \to \infty} A_{\tau}^{n} \pi_{\tau}
 \tau^{*~n}
\left(
\sum_{j \in {\cal A}} \sum_{w_k^{(j)}=i^{-1}}
-({\hako e}_j -A_{\tau}^{-1}{\bf f}(P_k^{(j)}w_k^{(j)}),j^*)
\right)\\
&= & 
 \cup_{j \in {\cal A}} \cup_{w_k^{(j)}=i^{-1}} 
(-A_{\tau}^{-1} \pi_{\tau}{\bf f}(P_k^{(j)}w_k^{(j)})+ {X}_{\tau}^{(j)} ),
\end{eqnarray*}
$\tau=\tau_2,\tau_4,~i \in {\cal A}$.
The other set equation for ${X'}_{\tau}^{(i)}$ is shown analogously.
\hfill$\Box$ \\

Secondly, to construct Rauzy fractals related to automorphism 
$\sigma=\sigma_1, \sigma_2, \sigma_4$, 
set seeds $\overline{\cal U}$ and $\overline{\cal U}'$ as
\begin{eqnarray*}
\overline{\cal U} &:=&
\left\{
\begin{array}{ll}
({\hako e}_1, 1^*) + ({\hako e}_2, 2^*) & \mbox{ if } \sigma=\sigma_2\\
({\hako o}, 1^*) + ({\hako e}_2, 2^*) & \mbox{ if } \sigma=\sigma_1,\sigma_4,\\
\end{array}
\right. \\
\overline{\cal U}' &:=&
\left\{
\begin{array}{ll}
({\hako o}, 1^*) + ({\hako o}, 2^*) & \mbox{ if } \sigma=\sigma_2\\
({\hako e}_1, 1^*) + ({\hako o}, 2^*) & \mbox{ if } \sigma=\sigma_1,\sigma_4.\\
\end{array}
\right. 
\end{eqnarray*}
Rauzy fractals related to $\sigma_t,~t=1,2,4$ are defined as follows.
\begin{definition}
\label{def:3-1}
The following limit sets exist in the sense of Lebesgue measure.
For $\sigma=\sigma_1,\sigma_2,\sigma_4$,
\begin{eqnarray*}
X_{\sigma}&:=& \lim_{n \to \infty} A_{\sigma}^n \pi_{\sigma}
 \sigma^{*~n}(\overline{\cal U}),\\
&:=& \lim_{n \to \infty} A_{\sigma}^n \pi_{\sigma}
 \sigma^{*~n}(\overline{\cal U}').
\end{eqnarray*}

In the case of (ii), for $i \in {\cal A}$,
\begin{eqnarray*}
X_{\sigma_2}^{(i)}&:=& \lim_{n \to \infty} A_{\sigma_2}^n \pi_{\sigma_2}
 \sigma_2^{*~n}({\hako e}_i, i^*),\\
{X'}_{\sigma_2}^{(i)}  &:=& \lim_{n \to \infty} A_{\sigma_2}^n \pi_{\sigma_2}
 \sigma_2^{*~n}({\hako o}, i^*).
\end{eqnarray*}

In the case of (i) and (iv), for $\sigma=\sigma_1,\sigma_4 $,
\begin{eqnarray*}
X_{\sigma}^{(1^{-1})}&:=& \lim_{n \to \infty} A_{\sigma}^n \pi_{\sigma}
 \sigma^{*~n}({\hako o}, 1^*),\\
X_{\sigma}^{(2)}&:=& \lim_{n \to \infty} A_{\sigma}^n \pi_{\sigma}
 \sigma^{*~n}({\hako e}_2, 2^*),\\
{X'}_{\sigma}^{(1^{-1})}  &:=& \lim_{n \to \infty} A_{\sigma}^n \pi_{\sigma}
 \sigma^{*~n}({\hako e}_1, 1^*),\\
{X'}_{\sigma}^{(2)}  &:=& \lim_{n \to \infty} A_{\sigma}^n \pi_{\sigma}
 \sigma^{*~n}({\hako o}, 2^*).
\end{eqnarray*}
\end{definition}

For each automorphism $\sigma_t,~t=1,2,4$, we set 
\begin{eqnarray*}
 \epsilon_1&:=&\left\{ 
\begin{array}{ll}
1 & \mbox{ if } \sigma=\sigma_2\\
-1 & \mbox{ if } \sigma=\sigma_1,\sigma_4,
\end{array}
\right.,\\
 \epsilon_2 &:= & 1,
\end{eqnarray*}
then 
$\delta^{-1}_t(i) \in \{ 1^{\epsilon_1}, 2^{\epsilon_2}\}^*$ for any $t=1,2,4,~i \in {\cal A}$.
Then we have the theorem which gives the relation between
Rauzy fractals $X_{\tau_t}^{(i)}$
and $X_{\sigma_t}^{(i^{\epsilon_i})}$, $t=1,2,4$.

\begin{theorem}
\label{theorem:3-1}
For $t=1,2,4$ and $i \in {\cal A}$, the following equations hold:
\begin{eqnarray*}
 X_{\sigma_t}^{(i^{\epsilon_i})}&=&
\cup_{j \in {\cal A}} \cup_{w_k^{(j)}=i^{\epsilon_i}}
(-\pi_{\sigma_t}{\bf f} (P_k^{(j)})+
A_{\delta_t}^{-1}X_{\tau_t}^{(j)}),\\
 {X'}_{\sigma_t}^{(i^{\epsilon_i})}&=&
\cup_{j \in {\cal A}} \cup_{w_k^{(j)}=i^{\epsilon_i}}
(\pi_{\sigma_t}{\bf f} (S_k^{(j)})+
A_{\delta_t}^{-1}{X'}_{\tau_t}^{(j)}),
\end{eqnarray*}
 where $\delta_t^{-1}(i)$ is written as $\delta_t^{-1}(i)=w_1^{(i)}\cdots w_k^{(i)}
\cdots w_{l^{(i)}}^{(i)}=P_k^{(i)} w_k^{(i)}S_k^{(i)} $.
The unions are disjoint in the sense of Lebesgue
 measure. Moreover,  
$ X_{\sigma_t}^{(i^{\epsilon_i})}$, 
${X'}_{\sigma_t}^{(i^{\epsilon_i})}$,
$ X_{\sigma_t}$, 
${X'}_{\sigma_t}$ are interval.
\end{theorem}

\begin{figure}[hbtp]
\setlength{\unitlength}{1mm}
The case (i):
\begin{center}
\epsfxsize=6cm
\epsfbox{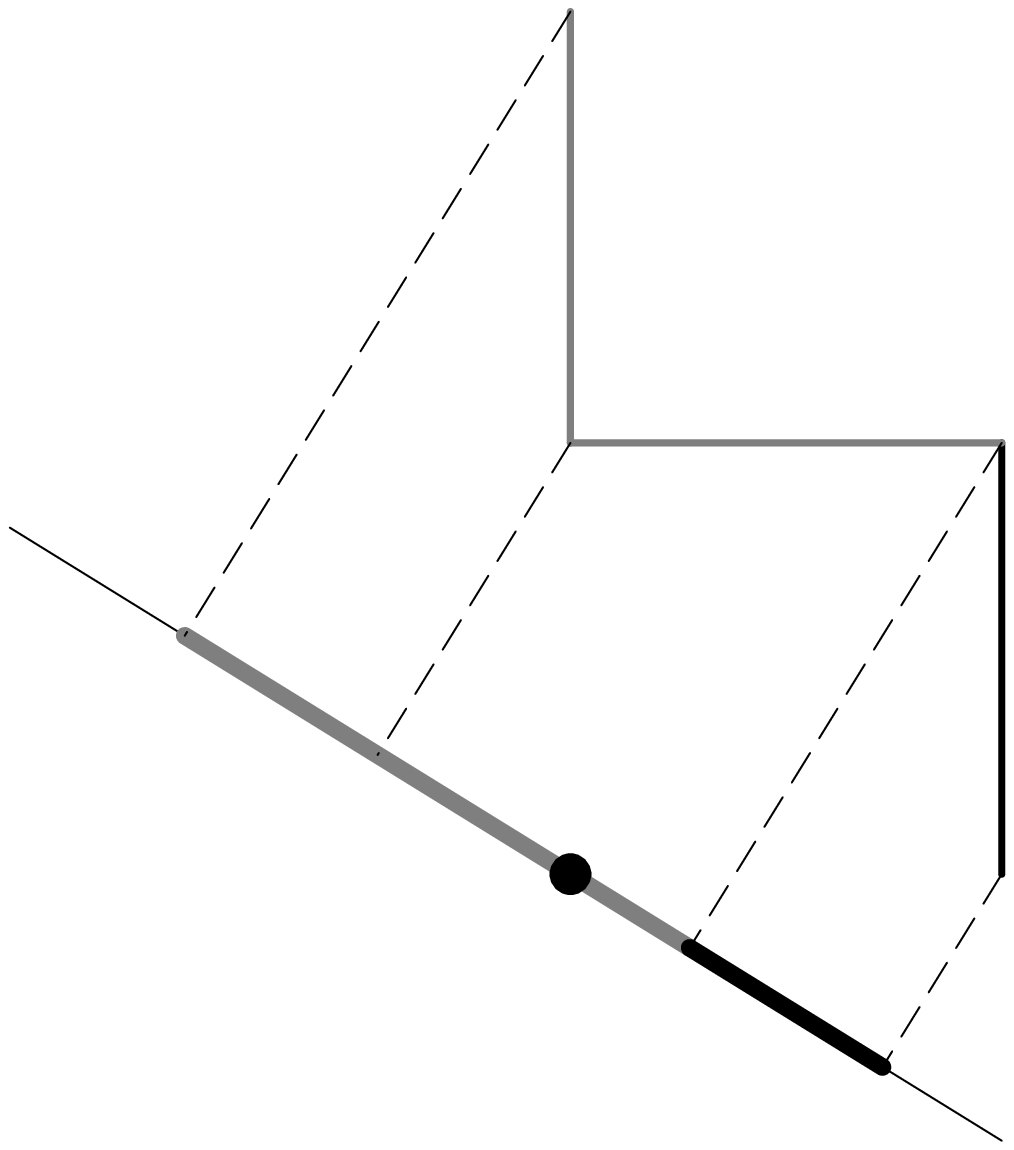}
~~~~
\epsfxsize=6cm
\epsfbox{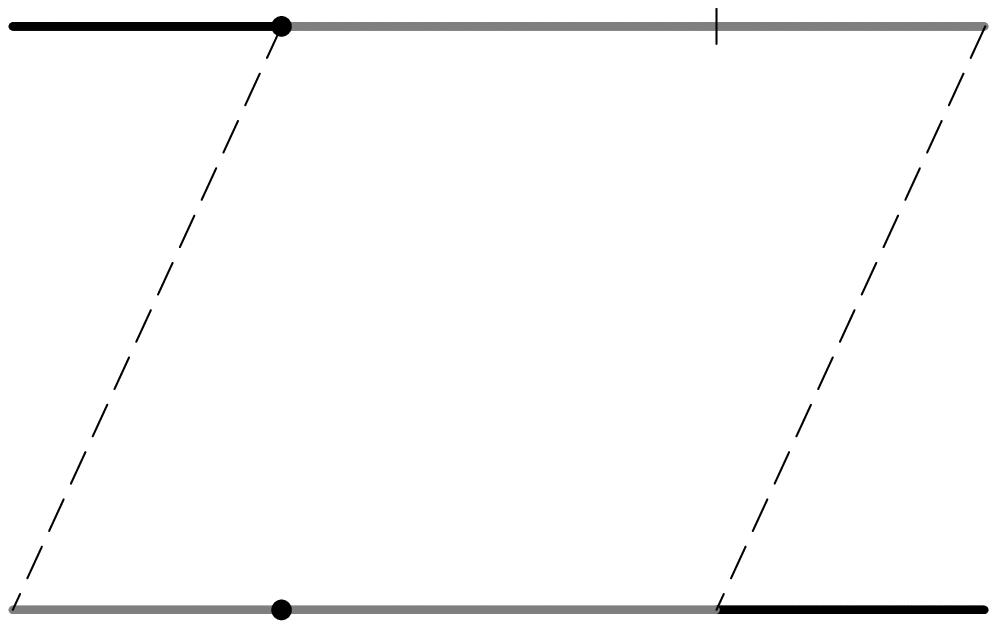}
\end{center}
\setlength{\unitlength}{1mm}
\begin{picture}(0,0)(0,0)
\put(10,50){$P_{\tau_1}$}
\put(0,36){$X_{\tau_1}^{(1)}$}
\put(0,30){$+\pi_{\tau_1}({\hako e}_2-{\hako
 e}_1)-{\hako h}$}
\put(25,22){$X_{\tau_1}^{(2)}-{\hako h}$}
\put(40,15){$X_{\tau_1}^{(1)}-{\hako h}$}
\put(85,50){$A_{\delta_1}^{-1}X_{\tau_1}^{(1)}$}
\put(105,50){$A_{\delta_1}^{-1}X_{\tau_1}^{(2)}$}
\put(125,55){$A_{\delta_1}^{-1}(X_{\tau_1}^{(1)}$}
\put(125,48){$+\pi_{\tau_1}({\hako e}_2-{\hako
 e}_1))$}
\put(85,37){$X_{\sigma_1}^{(1^{-1})}$}
\put(110,37){$X_{\sigma_1}^{(2)}$}
\end{picture}
The case (ii):
\begin{center}
\epsfxsize=6cm
\epsfbox{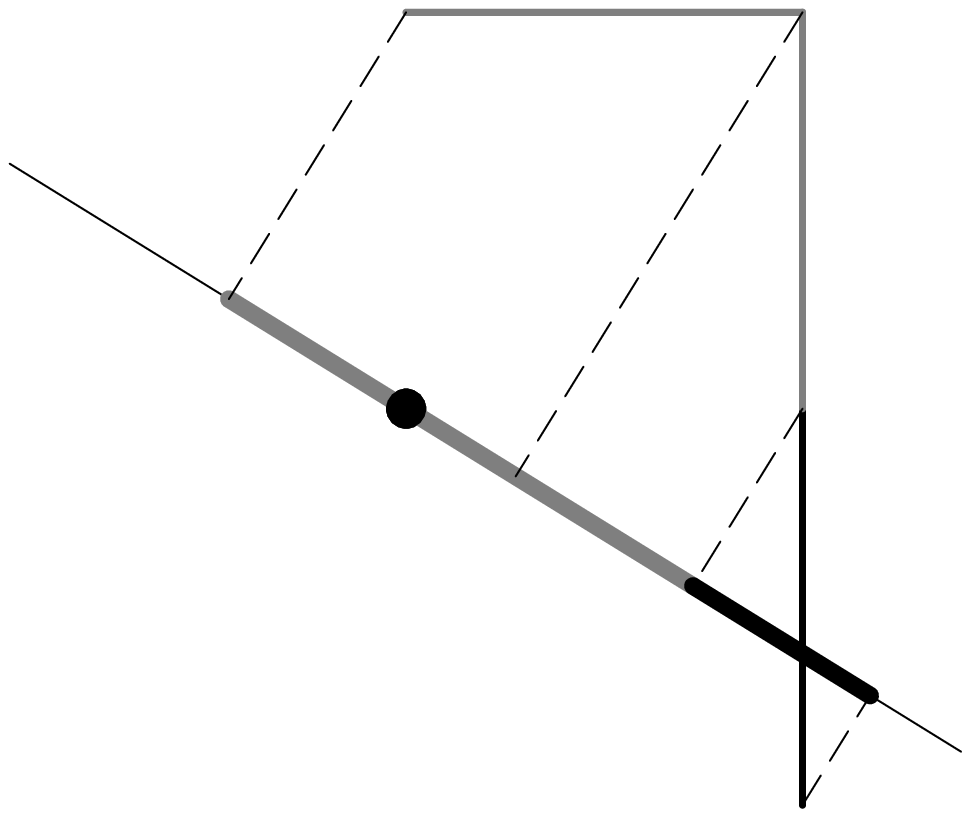}
~~~~
\epsfxsize=6cm
\epsfbox{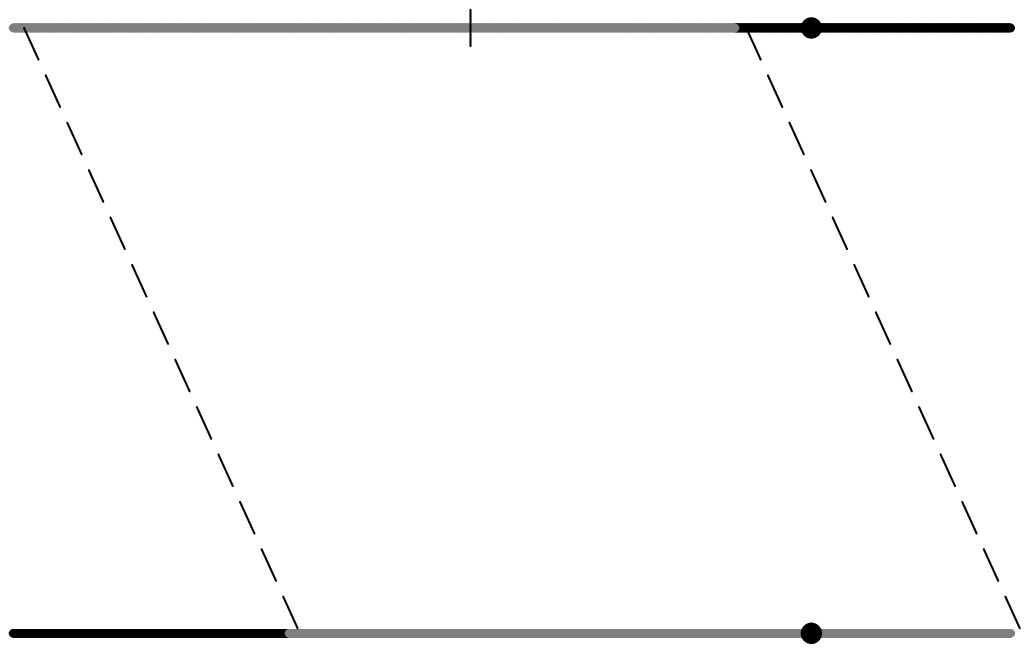}
\end{center}
\setlength{\unitlength}{1mm}
\begin{picture}(0,0)(0,0)
\put(10,50){$P_{\tau_2}$}
\put(15,30){$X_{\tau_2}^{(2)}-{\hako h}$}
\put(25,22){$X_{\tau_2}^{(1)}-{\hako h}$}
\put(30,15){$X_{\tau_2}^{(1)}-\pi_{\tau_2}{\hako e}_2-{\hako h}$}
\put(90,50){$A_{\delta_2}^{-1}X_{\tau_2}^{(2)}$}
\put(110,50){$A_{\delta_2}^{-1}X_{\tau_2}^{(1)}$}
\put(130,55){$A_{\delta_2}^{-1}(X_{\tau_2}^{(1)}$}
\put(130,48){$-\pi_{\tau_2}{\hako e}_2)$}
\put(105,37){$X_{\sigma_2}^{(2)}$}
\put(133,37){$X_{\sigma_2}^{(1)}$}
\end{picture}
The case (iv):
\begin{center}
\epsfxsize=6cm
\epsfbox{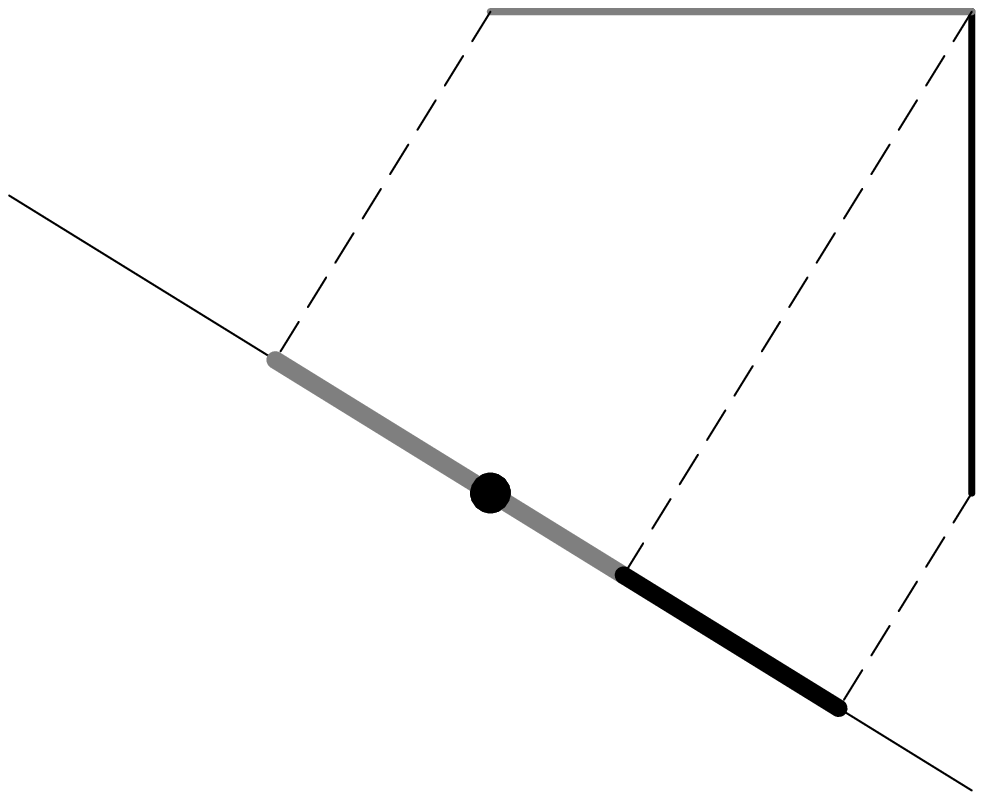}
~~~~
\epsfxsize=6cm
\epsfbox{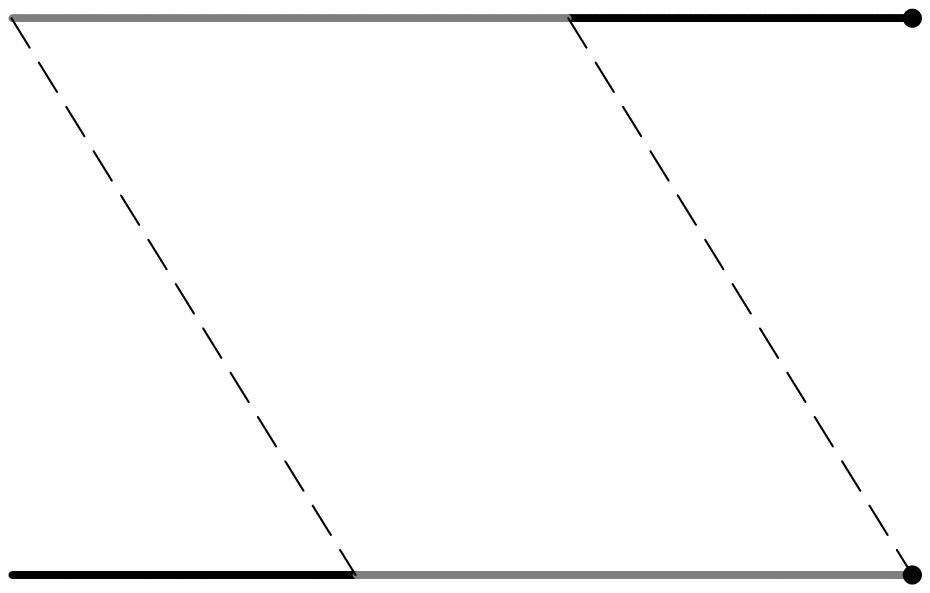}
\end{center}
\setlength{\unitlength}{1mm}
\begin{picture}(0,0)(0,0)
\put(10,50){$P_{\tau_4}$}
\put(25,22){$X_{\tau_4}^{(2)}-{\hako h}$}
\put(40,15){$X_{\tau_4}^{(1)}-{\hako h}$}
\put(95,50){$A_{\delta_4}^{-1}X_{\tau_4}^{(2)}$}
\put(127,50){$A_{\delta_4}^{-1}X_{\tau_4}^{(1)}$}
\put(103,37){$X_{\sigma_4}^{(2)}$}
\put(128,37){$X_{\sigma_4}^{(1^{-1})}$}
\end{picture}
\caption{$(\delta_t^{-1})^*(\overline{\cal U})$ and domain exchange
 transformations $T_{\sigma_t}$ on Rauzy fractals 
$X_{\sigma_t}^{(i^{\epsilon_i})},~i \in {\cal A},~t=1,2,4$}\label{fig:6}
\end{figure}

Proof.
Let us show the first equation for $\sigma=\sigma_2=\delta_2^{-1}
\circ \tau_2 \circ \delta_2$ in the case of
(ii). The other cases can be proved analogously.
By the definition of $ X_{\sigma}^{(i)}$, 
\begin{eqnarray*}
 X_{\sigma}^{(i)}
&=& \lim_{n \to \infty} A_{\sigma}^n \pi_{\sigma} \sigma^{*~n} 
({\hako e}_i,i^*)\\
&=& \lim_{n \to \infty} A_{\sigma}^n \pi_{\sigma} \delta^{*} \circ
\tau^{*~n}  \circ (\delta^{-1})^*
({\hako e}_i,i^*)\\
&=& \lim_{n \to \infty} A_{\sigma}^n \pi_{\sigma} \delta^* \circ
\tau^{*~n} \left\{ \sum_{j \in {\cal A}} \sum_{w_k^{(j)}=i} 
(A_{\delta}({\hako e}_i+{\bf f}(S_k^{(j)})),j^*) \right\}\\
&=&  \cup_{j \in {\cal A}} \cup_{w_k^{(j)}=i} 
\lim_{n \to \infty} A_{\sigma}^n \pi_{\sigma} \delta^* \circ
\tau^{*~n}
({\hako e}_j-A_{\delta}{\bf f}(P_k^{(j)}),j^*).
\end{eqnarray*}
Put
\[
 c_0:= \max_{i \in {\cal A}}
d_H \left(\pi_{\sigma} \delta^* ({\hako o},i^*),
\pi_{\sigma} A_\delta^{-1} ({\hako o},i^*) \right),
\]
where $d_H$ is the Hausdorff metric. From the property of the 
Hausdorff metric,
\begin{eqnarray*}
d_H \left(\pi_{\sigma} \delta^* \tau^{*~n} ({\hako o},i^*),
\pi_{\sigma} A_\delta^{-1} \tau^{*~n}({\hako o},i^*) \right)
& \le& c_0,\\
d_H \left(A_{\sigma}^n \pi_{\sigma} \delta^* \tau^{*~n} 
({\hako o},i^*),
A_{\sigma}^n \pi_{\sigma} A_\delta^{-1} \tau^{*~n}({\hako o},i^*)
    \right) 
&\le & c_0 |\lambda_{\sigma}'|^n,
\end{eqnarray*}
where $\lambda_{\sigma} '$ is the eigenvalue of $A_{\sigma}$
with $|\lambda_{\sigma}'|<1$. So 
\[
\lim_{n \to \infty}d_H \left(A_{\sigma}^n \pi_{\sigma} \delta^* \tau^{*~n}
({\hako o},i^*),
A_{\sigma}^n \pi_{\sigma} A_\delta^{-1} \tau^{*~n}({\hako o},i^*)
    \right) 
=0.
\]
By noticing the equality,
\[
A_{\delta}^{-1} \pi_{\tau} {\hako x} = \pi_{\sigma}A_{\delta}^{-1}
{\hako x},~
{\hako x} \in {\mbox{\boldmath $R$}}^2
\]
if $\sigma=\delta^{-1} \circ \tau \circ \delta$,
we have
\begin{eqnarray*}
\lim_{n \to \infty} A_{\sigma}^n \pi_{\sigma} \delta^* \circ
\tau^{*~n}({\hako e}_i-{\hako x},i^*)
&=& 
\lim_{n \to \infty} 
A_{\sigma}^n \pi_{\sigma} A_\delta^{-1} \tau^{*~n}
({\hako e}_i-{\hako x},i^*)\\
&=&
\lim_{n \to \infty} 
A_{\delta}^{-1} A_{\tau}^n \pi_{\tau} \tau^{*~n}
({\hako e}_i-{\hako x},i^*)\\
&=&
- \pi_{\sigma}A_{\delta}^{-1}{\hako x}
+ A_{\delta}^{-1} X_{\tau}^{(i)}.
\end{eqnarray*}
Therefore we have the set equation
\[
 X_{\sigma}^{(i)}
= \cup_{j \in {\cal A}} \cup_{w_k^{(j)}=i}
(-\pi_{\sigma}{\bf f} (P_k^{(j)})+
A_{\delta}^{-1}X_{\tau}^{(j)}).
\]
Next we show it is disjoint union and interval
through $A_{\delta} X_{\sigma}^{(i)}$. 
From Theorem \ref{theorem:0-2},
\begin{eqnarray*}
A_{\delta} X_{\sigma}^{(i)} 
&=& \cup_{j \in {\cal A}} \cup_{w_k^{(j)}=i}
(-\pi_{\tau}A_{\delta}{\bf f} (P_k^{(j)})+ X_{\tau}^{(j)})\\
&=& \cup_{j \in {\cal A}} \cup_{w_k^{(j)}=i}
\pi_{\tau}(- A_{\delta}{\bf f}(P_k^{(j)})+ {\hako e}_j,j^* )+{\hako h}\\
&=& \pi_{\tau}(\delta^{-1})^*({\hako e}_i,i^*)  +{\hako h}
\end{eqnarray*}
for some ${\hako h} \in P_{\tau}$.
It means we can obtain $A_{\delta} X_{\sigma}^{(i)}$ after
projection of $(\delta^{-1})^*({\hako e}_i,i^*)$ by $\pi_{\tau}$ and
translation by ${\hako h}$.
From the figure of $(\delta^{-1})^*({\hako e}_i,i^*) $, the union in
the equation are disjoint and $A_{\delta} X_{\sigma}^{(i)} $ is an
interval (see Figure \ref{fig:6}).
The other equation is shown analogously.
\hfill$\Box$\\

\begin{remark}
\label{remark:3-2}
From Figure \ref{fig:6}, $(\delta_t^{-1})^*(\overline{\cal U}) \notin
 {\cal G}_\tau^*$ in the case of (i), (ii) for $t=1,2$. Therefore
$\sigma^{*}(\overline{\cal U})%,~n \in {\mbox{\boldmath $N$}}
$ is not included in the stepped surface ${\cal S}_{\sigma_t}$. 
That is the reason why
we take the different seeds to construct the stepped surface.
\end{remark}
\begin{definition}
The domain exchange transformations $T_{\tau}$ on $X_{\tau}$ for
$\tau=\tau_1,\tau_2,\tau_4$
and $T_{\sigma}$ on $X_{\sigma}$ for $\sigma=\sigma_1,\sigma_2,\sigma_4$
are defined by
\[
\begin{array}{l}
T_{\tau}: X_{\tau} \to X_{\tau}\\
T_{\tau}({\hako x})={\hako x}- \pi_{\tau} {\bf f}(i) {\mbox ~if~ }
 {\hako x} \in  X_{\tau}^{(i)},
\end{array}
\]
and
\[
\begin{array}{l}
T_{\sigma}: X_{\sigma} \to X_{\sigma}\\
T_{\sigma}({\hako x})={\hako x}- \pi_{\sigma} {\bf f}(i^{\epsilon_i})
{\mbox ~if~ } {\hako x} \in  X_{\sigma}^{(i^{\epsilon_i})}.
\end{array}
\]
\end{definition}
By the definitions of $X_{\tau}^{(i)},~\tau=\tau_1,\tau_2,\tau_4,$
and $X_{\sigma}^{(i^{\epsilon_i})},~\sigma=\sigma_1,\sigma_2,\sigma_4$,
\[
X_{\tau}^{(i)}-\pi_{\tau}{\bf f}(i)={X'}_{\tau}^{(i)},~
X_{\sigma}^{(i^{\epsilon_i})}-\pi_{\sigma}{\bf f}(i^{\epsilon_i})
={X'}_{\sigma}^{(i{\epsilon_i})},
\]
therefore the domain exchange transformations are well-defined (see
Figure \ref{fig:6}).
\begin{theorem}
\label{theorem:3-2}
The measurable dynamical system $(X_{\sigma},T_{\sigma},\mu)$  with
Lebesgue measure $\mu$
has $\delta^{-1}$-structure with
 respect to the pair of partitions 
$\{ X_{\sigma}^{(i^{\epsilon_i})}|~i \in {\cal A} \}$,
$\{A_{\delta}^{-1} X_{\tau}^{(i)}|~i \in {\cal A}\}$ for
$\sigma=\sigma_t,~\tau=\tau_t,~\delta=\delta_t,~t=1,2,4$.
Moreover, $(X_{\sigma_1},T_{\sigma_1},\mu)$
 (resp. $(X_{\sigma_t},T_{\sigma_t},\mu),~t=2,4$)
has $\delta_1^{-1} \tau_1^n$-structure 
(resp. $\delta_t^{-1} \tau_t^{2n}$-structure)
with respect to the pair of partitions 
$\{X_{\sigma_1}^{(i^{\epsilon_i})}|~i \in {\cal A}\}$,
$\{A_{\delta_1}^{-1} A_{\tau_1}^n X_{\tau_1}^{(i)}|~i \in {\cal A} \}$
(resp. $\{ X_{\sigma_t}^{(i^{\epsilon_i})} |~i \in {\cal A} \}$,
$\{A_{\delta_t}^{-1} A_{\tau_t}^{2n} X_{\tau_t}^{(i)} |~i \in {\cal A}
 \}$) for any positive integer $n$ (see Figure \ref{fig:6}).
\end{theorem}

Proof.
From the set equation
\[
 X_{\sigma}^{(i^{\epsilon_i})}=
\cup_{j \in {\cal A}} \cup_{w_k^{(j)}=i^{\epsilon_i}}
(-\pi_{\sigma}{\bf f} (P_k^{(j)})+
A_{\delta}^{-1}X_{\tau}^{(j)}),
\] 
where $\delta^{-1}(i)$ is written as 
$\delta^{-1}(i)=w_1^{(i)}\cdots w_k^{(i)}
\ldots w_{l^{(i)}}^{(i)}=P_k^{(i)} w_k^{(i)}S_k^{(i)} $,
\begin{center}
\(
\begin{array}{l}
A_{\delta}^{-1}X_{\tau}^{(i)} \subset X_{\sigma}^{(w_1^{(i)})},\\
A_{\delta}^{-1}X_{\tau}^{(i)} - \pi_{\sigma} {\bf f}(w_1^{(i)})
\subset X_{\sigma}^{(w_2^{(i)})},\\
\ldots\\
A_{\delta}^{-1}X_{\tau}^{(i)} - \pi_{\sigma} {\bf f}
(w_1^{(i)}\cdots w_{l^{(i)}-1}^{(i)})
\subset X_{\sigma}^{(w_{l^{(i)}}^{(i)})},\\
A_{\delta}^{-1}X_{\tau}^{(i)} - \pi_{\sigma} A_{\delta}^{-1}{\bf f}(i)
= A_{\delta}^{-1}(X_{\tau}^{(i)} - \pi_{\tau} {\bf f}(i))
= A_{\delta}^{-1}{X'}_{\tau}^{(i)}.
\end{array}
\)
\end{center}
So
$(X_{\sigma},T_{\sigma},\mu)$ has $\delta^{-1}$-structure.
From
\[
 T_{\sigma}|_{A_{\delta}^{-1} X_{\tau}}
({\hako x}) =A_{\delta}^{-1} \circ T_{\tau} \circ A_{\delta} ({\hako
x}),~{\hako x} \in X_{\sigma},
\]
the induced transformation $ T_{\sigma}|_{A_{\delta}^{-1} X_{\tau}}$ is 
conjugate to $ T_{\tau} $.
Recall that $T_{\tau_1}$ (resp. $T_{\tau_t},~t=2,4$) has $\tau_1$-structure 
 (resp. $\tau_t^2$-structure ) with respect to the pair of partitions
$\{X_{\tau_1}^{(i)}|~i \in {\cal A} \}$,
$\{A_{\tau_1}^n X_{\tau_1}^{(i)}|~i \in {\cal A} \}$
(resp. $\{X_{\tau_t}^{(i)}|~i \in {\cal A} \}$,
$\{A_{\tau_t}^{2n} X_{\tau_1}^{(i)}|~i \in {\cal A} \}$) by Theorem 
\ref{theorem:0-1}.
Thus the last part is proved.
\hfill$\Box$\\

One-sided sequence $\omega$ is called a {\it fixed point} for $\sigma$
if $\sigma (\omega)=\omega$; and $\omega$ is called a {\it periodic point} 
with period $n$ if $\sigma^n (\omega)=\omega$.
The substitution $\tau_1$ has a fixed point and the alternative
substitutions $\tau_t,~t=2,4$ have periodic points of period 2.
We denote the fixed point $\lim_{n \to \infty}\tau_1^n(2)$  by
$\omega_{\tau_1}$, and the periodic points 
$\lim_{n \to \infty}\tau_t^{2n}(2),~t=2,4$ by
$\omega_{\tau_t}$, where these limits exist in the sense of the product topology.
Let us define the one-sided sequences 
\[
\omega_{\sigma_1}  := \lim_{n \to \infty}\sigma_1^n(2),~
\omega_{\sigma_2}  := \lim_{n \to \infty}\sigma_2^{2n}(2),~
\omega_{\sigma_4}  := \lim_{n \to \infty}\sigma_4^{2n}(2).
\]
These fixed point or periodic points $\omega_{\sigma_t},~t=1,2,4$ are given by
\[
\omega_{\sigma_t}  =
 \delta_t^{-1}
( \omega_{\tau_t}) \in \{1^{\epsilon_1},2^{\epsilon_2}\}^{\mbox{\boldmath $N$}}.
\]
The one-sided sequences $\omega_{\sigma_t},\omega_{\tau_t},~t=1,2,4$
are written as
\begin{eqnarray*}
\omega_{\sigma_t}&=&s_0 s_1 \cdots s_k \cdots,\\
\omega_{\tau_t}&=&t_0 t_1 \cdots t_k \cdots. 
\end{eqnarray*}
%and $\omega_{\sigma_1}$ is a fixed point for $\sigma_1$ and
%$\omega_{\sigma_i}, i=2,4$ are periodic points of period 2 for
%$\sigma_i$. 
Since 
$\tau_t^{*~2n}({\hako e}_2, 2^*)$ includes $({\hako e}_2, 2^*)$,
$t=1,2,4$ for any positive integer $n$, and
the origin point ${\hako o} \in \pi_{\tau_t}({\hako e}_2, 2^*)$, 
so ${\hako o} \in X_{\tau_t}^{(2)}$.
The orbit of the origin point by $T_{\sigma_t}$ is described by 
a fixed point or a periodic point of $\sigma_t$ by Theorem \ref{theorem:3-2}.
\begin{corollary}
For $\sigma=\sigma_1,\sigma_2,\sigma_4$, 
\[
T_{\sigma}^k({\hako o}) \in  X_{\sigma}^{(s_k)},~k=0,1,\cdots. 
\]
\end{corollary}

%%%%%%%%%%%%%%%%%%
Finally we will see that Rauzy fractals related to
$\sigma_t,~t=1,2,4$ are also given by the fixed point or the
periodic point $\omega_{\sigma_t}$  as we saw in Section 0.

For a substitution or an alternative substitution $\tau=\tau_t,~t=1,2,4$,
put
\begin{eqnarray*}
Y_{\tau} &:=& \{ -\pi_{\tau} {\bf f} (t_0 t_1 \cdots t_k)|~k \ge 0 \},\\
Y_{\tau}^{(i)} &:=& \{- \pi_{\tau} {\bf f} (t_0 t_1 \cdots t_{k-1})|~
k \ge 0,~t_k=i \},\\
Y'{}_{\tau}^{(i)} &:=& \{- \pi_{\tau} {\bf f} (t_0 t_1 \cdots t_{k})|~
k \ge 0,~t_k=i \}.
\end{eqnarray*}
Since $\tau_1, \tau_2^2,\tau_4^2$ are substitutions and the equality
(\ref{atomic1}), (\ref{atomic2}) in Section 0 and Remark \ref{remark:3-1}, we have
\begin{eqnarray*}
 X_{\tau}&=&\overline {Y_{\tau}},\\
 X_{\tau}^{(i)}&=& X_{\tau^2}^{(i)}=\overline {Y_{\tau}^{(i)}}, \\
{X'}_{\tau}^{(i)}&=&{X'}_{\tau^2}^{(i)}=\overline{ {Y'}_{\tau}^{(i)}}.
\end{eqnarray*}
%where ${\overline S}$ means a closure of the set $S$.

%\begin{proposition}(?????)
%A substitution $\tau$ is a primitive invertible substitution if and
%only if atomic surfaces  $X_{\tau},X_{\tau}^{(i)}, X'{}_{\tau}^{(i)}$ of  $\tau$
%are interval, where a substitution is refer to be primitive if there
%exists $n$ such that $A_\tau^n >0$.
%\end{proposition}
%Since $\tau=\tau_1,\tau_2^2,\tau_4^2$ are primitive invertible
%substitutions, atomic surfaces  $X_{\tau},X_{\tau}^{(i)},
%X'{}_{\tau}^{(i)}$ are interval from the Proposition ????.\\
%
%For $\sigma =\sigma_i, i=1,2,4$, we define atomic surfaces in the same
%way.
%
For automorphisms $\sigma=\sigma_1,\sigma_2, \sigma_4$, we have the same
result.
\begin{theorem}
\label{theorem:3-3}
%The fixed point or the periodic points $\omega_\sigma$
%of $\sigma$ are written as $\omega_{\sigma}=s_0 s_1 \cdots s_k \cdots$
%for $\sigma =\sigma_t, t=1,2,4$, and put 
For $\sigma=\sigma_1,\sigma_2, \sigma_4$ and $i \in {\cal A}$, put
\begin{eqnarray*}
Y_{\sigma} &:=& \{- \pi_{\sigma} {\bf f} (s_0 s_1 \cdots s_k)|~k \ge 0 \}\\
Y_{\sigma}^{(i^{\epsilon_i})} &:=& 
\{- \pi_{\sigma} {\bf f} (s_0 s_1 \cdots s_{k-1})|
~k \ge 0,~s_k=i^{\epsilon_i} \}\\
Y'{}_{\sigma}^{(i^{\epsilon_i})} &:=& 
\{- \pi_{\sigma} {\bf f} (s_0 s_1 \cdots s_{k})|~k \ge 0,~s_k=i^{\epsilon_i} \}.
\end{eqnarray*}
Then the following equalities hold:
\[
 X_{\sigma}=\overline{ Y_{\sigma}},~ X_{\sigma}^{(i^{\epsilon_i})}=\overline{
 Y_{\sigma}^{(i^{\epsilon_i})}},~ 
X'{}_{\sigma}^{(i^{\epsilon_i})}=\overline{
 Y'{}_{\sigma}^{(i^{\epsilon_i})}},~
~i \in {\cal A}.
\]
\end{theorem}

Proof.
By Theorem \ref{theorem:3-1}, 
to prove the equality $ X_{\sigma}^{(i^{\epsilon_i})}=\overline{
 Y_{\sigma}^{(i^{\epsilon_i})}}$,
 it is enough to show that
\[
  Y_{\sigma}^{(i^{\epsilon_i})}=
\cup_{j \in {\cal A}} \cup_{w_k^{(j)}=i^{\epsilon_i}}
(-\pi_{\sigma}{\bf f} (P_k^{(j)})+
A_{\delta}^{-1}Y_{\tau}^{(j)}),
\]
where $\delta^{-1}(i)=w_1^{(i)}\cdots w_k^{(i)} \cdots w_{l^{(i)}}^{(i)}
 =P_k^{(i)} w_k^{(i)}S_k^{(i)},~i \in {\cal A}$.
Take $-\pi_{\sigma}{\bf f}(s_0 s_1 \cdots s_{k-1}) \in
Y_{\sigma}^{(i^{\epsilon_i})}$ such that $ s_k=i^{\epsilon_i}$.
There exist $k_1, k_2$ such that
\[
s_0 s_1 \cdots s_{k-1}=\delta^{-1}(t_0 t_1 \cdots 
t_{k_1 -1})P_{k_2}^{(t_{k_1})},~
w_{k_2}^{(t_{k_1})}=i^{\epsilon_i}.
\]
Then
\begin{eqnarray*}
-\pi_\sigma {\bf f}(s_0 s_1 \cdots s_{k-1})&=&
-\pi_\sigma (A_{\delta}^{-1} {\bf f}(t_0 t_1 \cdots t_{k_1-1}) +
{\bf f}(P_{k_2}^{(t_{k_1})}))\\
&=&
-\pi_\sigma {\bf f}(P_{k_2}^{(t_{k_1})})
- A_{\delta}^{-1} \pi_\tau {\bf f}(t_0 t_1 \cdots t_{k_1-1}) ,
\end{eqnarray*}
and $-\pi_\sigma {\bf f}(s_0 s_1 \cdots s_{k-1}) \in 
- \pi_\sigma 
{\bf f}(P_{k_2}^{(t_{k_1})}) +A_{\delta}^{-1}Y_{\tau}^{(t_{k_1})} $.
Therefore we have
\[
  Y_{\sigma}^{(i^{\epsilon_i})} \subset
\cup_{j \in {\cal A}} \cup_{w_k^{(j)}=i^{\epsilon_i}}
(-\pi_{\sigma}{\bf f} (P_k^{(j)})+
A_{\delta}^{-1}Y_{\tau}^{(j)}).
\]
The opposite inclusive relation can be shown easily.
\hfill$\Box$